\documentclass[12pt,a4paper]{article}
\usepackage[utf8]{inputenc}
\usepackage[english]{babel}
\usepackage[ruled,section]{algorithm}
\usepackage{makeidx}
\usepackage[font=small]{caption}
\usepackage{amsmath}
\usepackage{amsthm}
\usepackage{amsfonts}
\usepackage{amssymb}
\usepackage{mathrsfs}
\usepackage{graphicx}
\usepackage[margin=25mm]{geometry}
\usepackage{enumerate}
\usepackage[plainpages=false,pagebackref=false]{hyperref}
\usepackage{color}
\usepackage{authblk}
\usepackage{hyperref}
\makeindex
\hypersetup{pdfpagelabels,hyperindex,colorlinks=true,breaklinks=true,bookmarks=true,linkcolor=black,citecolor=black,
urlcolor=black}
\newcommand{\real}{\mathbb{R}}
\newcommand{\n}{\mathbb{N}}
\newcommand{\rn}{ {\mathbb{R}^n} }
\renewcommand{\d}{\; \mathrm{d}}
\numberwithin{equation}{section}
\newtheorem{teo}{Theorem}[section]
\newtheorem{lem}[teo]{Lemma}
\newtheorem{rmk}[teo]{Remark}

\newtheorem{defin}[teo]{Definition}
\newtheorem{corol}[teo]{Corollary}
\newtheorem{prop}[teo]{Proposition}

\newtheorem{claim}[teo]{Claim}

\begin{document}


\title{\bf\Large A priori bounds and multiplicity for fully nonlinear equations with quadratic growth in the gradient}
\author[1]{Gabrielle Nornberg\footnote{gabrielle@mat.puc-rio.br, supported by Capes PROEX/PDSE grant 88881.134627/2016-01.}}
\author[1]{Boyan Sirakov\footnote{bsirakov@mat.puc-rio.br}}
\affil[1]{Pontifícia Universidade Católica do Rio de Janeiro, Brazil}
\date{}

\maketitle

{\small\noindent{\bf{Abstract.}} We consider fully nonlinear uniformly elliptic equations with quadratic growth in the gradient, such as
$$
-F(x,u,Du,D^2u) =\lambda c(x)u+\langle M(x)D u, D u \rangle +h(x)
$$
in a bounded domain with a Dirichlet boundary condition; here $\lambda \in\real$,  $c,\, h \in L^p(\Omega)$, $p>n\geq 1$, $c\gneqq 0$ and the matrix $M$ satisfies $0<\mu_1 I\leq M\leq \mu_2 I$. Recently this problem was studied in the ``coercive" case $\lambda c\le0$, where uniqueness of solutions can be expected; and it was conjectured  that the solution set is more complex for noncoercive equations. This conjecture was verified in~2015 by Arcoya, de Coster, Jeanjean and Tanaka for equations in divergence form, by exploiting the integral formulation of the problem. Here we show that similar phenomena occur for general, even fully nonlinear, equations in nondivergence form.  We use different techniques based on the maximum principle.

We develop a new method to obtain the crucial uniform a priori bounds, which permit to us to use degree theory. This method is based on basic regularity estimates such as half-Harnack inequalities, and on a Vázquez type strong maximum principle for our kind of equations.}

\section{Introduction}\label{Introduction}

This paper studies nonlinear uniformly elliptic problems of the following form
\begin{align}\label{Plambda} \tag*{$(P_\lambda)$}
\left\{
\begin{array}{rclcc}
-F(x,u,Du,D^2u) &=&\lambda c(x)u+\langle M(x)D u, D u \rangle +h(x) &\mbox{in} & \Omega  \\
u&=& 0 &\mbox{on} & \partial\Omega
\end{array}
\right.
\end{align}
where $\Omega$ is a bounded $C^{1,1}$ domain in $\rn$, $\lambda\in \mathbb{R}$, $n\geq 1$, $c,h\in L^p(\Omega)$, $M$ is a bounded matrix, and $F$ is a fully nonlinear uniformly elliptic operator of Isaacs type. A particular case, for which all our results are new as well, is when $F$ is a linear operator in nondivergence form i.e. $F(x,u,Du,D^2u)= a_{ij}(x)\partial_{ij}u + b_i(x)\partial_i u$.

A remarkable feature of this class of problems is that the second order terms and the gradient term have the same scaling with respect to dilations, so the second order term is not dominating when we zoom into a given point. This type of gradient dependence is usually named ``natural" in the literature, and is the object of extensive study. Another important property of \ref{Plambda} is the invariance of this class of equations with respect to diffeomorphic changes of the spatial variable $x$ and the dependent variable $u$.

The importance of these problems has long been recognized, at least since the classical works of Kazdan-Kramer \cite{KK78} and Boccardo-Murat-Puel \cite{BMP1}, \cite{BMP2}. The latter is a rather complete study of  solvability of strictly coercive equations in divergence form, that is, when $F$ is the divergence of an expression of $x$, $u$, and $Du$. Strictly coercive for \ref{Plambda} means that $\lambda c(x)<<0$, and then uniqueness of solutions is to be expected (see \cite{BBGK}, \cite{BM}). For weakly coercive equations (for instance, when $c\equiv0$), existence and uniqueness can be proved only under a smallness assumption on $c$ and $M$, as was first observed by Ferone-Murat \cite{FM}. All these works use the weak integral formulation of PDEs in divergence form.

The second author showed in \cite{arma2010} that the same type of existence and uniqueness results can be proved for general coercive (i.e. proper) equations in nondivergence form, by using techniques based on the maximum principle. In that paper it was also observed, for the first time and with a rather specific and simple example ($F=\Delta$, $c=1$, $M=I$, $h=0$), that the solution set may be very different in the nonproper case $c>0$, and in particular more than one solution may appear.

In the recent years appeared a series of papers which unveil the complex nature of the solution set for noncoercive equations, in the particular case when $F$ is the Laplacian. In \cite{JB} the authors used a mountain pass argument in the case $M=const. I$, when a classical exponential change reduces the equation to a semilinear one; later Arcoya et al~\cite{ACJT} developed a method based on degree theory which applies to general gradient terms, and this work was extended and completed by Jeanjean and de Coster \cite{CJ} who gave a description of the solution set in terms of the parameter $\lambda$. Souplet \cite{Souplet} showed that the study of \ref{Plambda} is even more difficult if $M$ is allowed to vanish somewhere (note that in any case a hypothesis which prevents $M\equiv0$ is necessary). In all these works the crucial a priori bounds for $u$ in the $L^\infty$-norm rely on the fact that the second order operator is the Laplacian (or a divergence form operator). A result with the $p$-Laplacian is found in \cite{CF}.

It is our goal here to perform a similar study for general operators in nondivergence form, and extend the results from \cite{arma2010} to noncoercive equations. We note that nondivergence (fully nonlinear) equations with natural growth are particularly relevant for applications, since problems with such growth in the gradient are  abundant in control and game theory, and more recently in mean-field problems, where Hamilton-Jacobi-Bellman and Isaacs operators appear  as infinitesimal generators of the underlying stochastic processes.

Our methods for obtaining a priori bounds in the uniform norm are (necessarily) very different from those in the preceding works, since for us no integral formulation of the equation is available. We prove that solutions of \ref{Plambda} are bounded from above by a new method, based on some standard estimates from regularity theory, such as half-Harnack inequalities, and their recent boundary extensions in \cite{B2016}. On the other hand, lower bounds are shown to be equivalent, somewhat surprisingly, to a Vázquez type strong maximum principle for our equations, which we also establish.

\vspace{0.1cm}
The paper is organized as follows. The next section contains the statements of our results. In the preliminary section \ref{Preliminaries} we recall some known results that will be used along the text. In section \ref{secao fixed point problem} we construct and study an auxiliary fixed point problem in order to obtain the existence statements in theorems \ref{th1.1,1.2,1.3}--\ref{th1.4}, via degree theory. The core of the paper is section \ref{secao a priori bounds}, where we prove a priori bounds in the uniform norm for the solutions of the noncoercive problem \ref{Plambda}.  Section \ref{secao proofs} is devoted to the proof of the main theorems.

\section{Main Results}

In this section we give the precise hypotheses and statements of our results.

In \ref{Plambda} we assume that
$c\gneqq 0$ -- so ``coercive" or ``proper" corresponds to $\lambda\le0$, and ``noncoercive" to $\lambda>0$. Our main goal is to give some description of the solution set of \ref{Plambda} when the parameter $\lambda$ is positive.

We assume that the matrix $M$ satisfies the nondegeneracy condition
\begin{align}\label{M}
\tag*{$(M)$}  \mu_1 I \leq M(x)\leq \mu_2 I\;\;\;\mathrm{a.e.\;\; in\;} \Omega
\end{align}
for some $\mu_1, \mu_2>0$, and that \ref{Plambda} has the following structure
\begin{align} \label{SC}
 \mathcal{M}^- (X-Y)-b(x)|\vec{p}-\vec{q}|-d(x)\omega((r-s)^+)
\leq F(x,r,\vec{p},X) - F(x,s,\vec{q},Y) \nonumber\\
\leq \mathcal{M}^+ (X-Y)+b(x)|\vec{p}-\vec{q}|+d(x)\omega((s-r)^+) \textrm{ a.e. } x\in \Omega \tag{SC} \\
 F(\cdot,0,0,0)\equiv 0\, , \;\;\; b,\, d,\, c,\,h \in L^p (\Omega), p>n, \,b,d\geq 0, \; \omega \textrm{ a modulus of continuity}. \nonumber
\end{align}
For most of our results we will need the following stronger assumption
\begin{align} \label{SC0}
\mathcal{M}^- (X-Y)-b(x)|\vec{p}-\vec{q}| \leq F(x,r,\vec{p},X) - F(x,s,\vec{q},Y) \nonumber\\
\leq \mathcal{M}^+ (X-Y)+b(x)|\vec{p}-\vec{q}| \;\textrm{ a.e. } x\in\Omega \tag*{$(\mathrm{SC})_0$}\\
 F(\cdot,0,0,0)\equiv 0 \, ,\qquad \;\; b,\,h \in L^\infty (\Omega) , \, b\geq 0.\nonumber
\end{align}
We note that a very particular case of the last hypothesis appears when $F$ is a general linear operator, but we can go much further, allowing $F$ to be an arbitrary supremum or infimum of such linear operators, i.e. a Hamilton-Jacobi-Bellman (HJB) operator, and even to be a sup-inf of linear operators (Isaacs operator). We denote with $\mathcal{M}^\pm$ the extremal Pucci operators (see the next section).

We will also assume that for some $\theta>0 \, , \;r_0>0$ and all $x_0\in \overline{\Omega}$
\begin{align}\label{Hbeta}\tag{$H_\beta$}\left( \frac{1}{r^n} \int_{B_r(x_0)\cap\Omega} \beta_F (x,x_0)^p \right)^{\frac{1}{p}} \leq \theta , \quad\mbox{ for all }\;r\leq r_0
\end{align}
where $\beta_F (x,x_0):= \sup\{{(\|X\|+1)}^{-1}|F(x,0,0,X)-F(x_0,0,0,X)|\}$, taken over all symmetric matrices $X$.
This  is satisfied, for instance, if $F(x,0,0,X)$ is continuous in $x\in \overline{\Omega}$ (if $F$ is linear this means $a_{ij}(x)$ are continuous). The conditions \ref{M}-\eqref{SC}-\eqref{Hbeta} guarantee that the viscosity solutions of \ref{Plambda} have global  $C^{1,\alpha}$-regularity and estimates. This was proved in \cite{regularidade}, building on and extending the previous works \cite{MilSilv}, \cite{Winter}, \cite{BS}.

Solutions of the Dirichlet problem \ref{Plambda} are understood in the $L^p$-viscosity sense (see the next section) and belong to $C(\overline{\Omega})$. Thus, we study and prove multiplicity of bounded solutions. We note that multiple unbounded solutions can easily be found  for simple equations with natural growth. For instance, in \cite{ADP2006} it was observed that $\Delta u=|D u|^2$ admits infinitely many weak solutions in $ W^{1,2}_0 (B_1)$, namely $u_k=\mathrm{ln}( ({|x|^{2-n}-k})({1-k})^{-1} )$, $0\leq k<1$, in the case $n>2$.

We recall that strong solutions are functions in $W^{2,p}_{\mathrm{loc}} (\Omega)$ which satisfy the equation almost everywhere. Strong solutions are viscosity solutions \cite{KSweakharnack}. Conversely, it is known that if $F$ is for instance convex in the matrix $X$ and satisfies \ref{SC0} (such are the HJB operators), then viscosity solutions are strong \cite{regularidade}, and the convexity assumption can be removed in some cases  but not in general -- see \cite{EE}. For some of our results we will need to assume that viscosity solutions of \ref{Plambda} are strong -- see \eqref{Hstrong} below.

Since we want to study the way the nature of the solution set changes when we go from negative to positive zero order term, we will naturally assume that the problem with $\lambda =0$ has a solution. We also assume that the Dirichlet problem for $F$ is uniquely solvable, so that we can concentrate on the way the coefficients $c$ and $M$ influence the solvability.
We now summarize these conditions on $F$. First, we assume that
\begin{align}\label{H0}
\tag*{$(H_0)$} \textrm{the problem (}P_0\textrm{) has a strong solution }u_0.
\end{align}
Further, setting $F[u]:=F(x,u,Du,D^2u)$, we assume that for each $f\in L^p(\Omega)$,
\begin{align}\label{ExistUnic M bem definido} \tag{$H_1$} \mathrm{there\;exists\;a\;unique}\,L^p\,\mathrm{viscosity\;solution\;of}\;
\begin{cases}
-F[u]=f(x) &\textrm{in} \;\; \Omega  \\
\hspace{0.9cm} u=0  &\textrm{on} \;\; \partial\Omega.
\end{cases}
\end{align}
Given $c,h$ for which we study \ref{Plambda}, if $(\overline{P}_\lambda )$ denotes the problem \ref{Plambda} with $c$ and $h$ replaced by $\overline{c}$ and $\overline{h}$, we sometimes require that  $L^p$-viscosity solutions $\overline{u}_\lambda$ of $(\overline{P}_\lambda )$ are such that
\begin{align} \label{Hstrong} \tag{$H_2$}
\overline{u}_\lambda\in W^{2,p} (\Omega), \; \mathrm{for\;every\;} 0 \leq\overline{c}\leq c\,,\;|\overline{h}|\leq |h|+1+c.
\end{align}

We observe that, by Theorem 1(iii) of \cite{arma2010}, the function $u_0$ is the unique $L^p$-viscosity solution of $(P_0)$. Theorem 1(ii) of \cite{arma2010} shows that \ref{H0} holds for instance if $Mh$ has small $L^p$-norm (examples showing that in general this hypothesis cannot be removed are also found in that paper). Recall also that  \eqref{ExistUnic M bem definido} and \eqref{Hstrong} are both true if $F$ satisfies \ref{SC0} and is convex or concave in $X$, by the results in  \cite{Winter} and  \cite{regularidade}.\medskip

We now state our  main results. The following theorem contains a crucial uniform estimate for solutions of \ref{Plambda}, which is both important in itself and instrumental for the existence statements below.

\begin{teo} \label{apriori}
Let $\Omega\in C^{1,1}$ be a bounded domain. Suppose \ref{SC0}, \ref{H0} hold and let $\Lambda_1, \, \Lambda_2$ with $0<\Lambda_1<\Lambda_2$. Then every $L^p$-viscosity solution  $u$ of \ref{Plambda} satisfies
$$
\|u^-\|_{L^\infty(\Omega)} \leq C\, , \;\textrm{ for all } \lambda\in [0,\Lambda_2],\qquad\|u^+\|_{L^\infty(\Omega)} \leq C\, , \;\textrm{ for all } \lambda\in [\Lambda_1,\Lambda_2],
$$
where $C$ depends on $n,p,\mu_1, \Omega, \Lambda_1,\Lambda_2,  \|b\|_{L^\infty(\Omega)},\|c\|_{L^\infty(\Omega)}, \|h\|_{L^\infty(\Omega)},\|u_0\|_{L^\infty(\Omega)},\lambda_P,\Lambda_P$, the $C^{1,1}$ diffeomorphism that describes the boundary, and  the set where $c>0$.
\end{teo}
It will be clear from the theorems below that the restrictions on $\lambda$ cannot be removed.

As in previous works, we use the following order in the space $E:=C^1(\overline{\Omega})$.

\begin{defin}\label{def1.2}
Let $u,v\in E$. We say that $u\ll v$ if for every $x \in \Omega$ we have $u(x)<v(x)$ and for $x_0\in\partial\Omega$ we have either $u(x_0)<v(x_0)$, or $u(x_0)=v(x_0)$ and $\partial_\nu u (x_0)<\partial_\nu v (x_0)$,  where $\vec{\nu}$ is the interior unit normal to $\partial \Omega$. \end{defin}

\begin{teo} \label{th1.1,1.2,1.3}
Assume \eqref{SC},  \eqref{Hbeta}, \ref{H0}, and \eqref{ExistUnic M bem definido}.

1. Then, for $\lambda\leq 0$, the problem \ref{Plambda} has an $L^p$-viscosity solution $u_\lambda$ that converges to $u_0$ in $E$ as $\lambda\rightarrow 0^-$. Moreover, the set
$$
\Sigma = \{ \,(\lambda , u) \in \real \times E\, ; \, u \;\, \textrm{solves \ref{Plambda}} \,\}
$$
possesses an unbounded component $\mathcal{C}^+\subset [0,+\infty]\times E$ such that $\mathcal{C}^+\cap ( \{0\}\times E )=\{u_0\}$.

From now on we assume \ref{SC0}.

2. The component from 1. is such that
\begin{enumerate}[(i)]
\item either it bifurcates from infinity to the right of the axis $\lambda =0$ with the corresponding solutions having a positive part blowing up to infinity in $C (\overline{\Omega})$ as $\lambda\rightarrow 0^+$;
\item or  its projection on the $\lambda$ axis is $[0,+\infty)$.
\end{enumerate}

3. There exists $\bar{\lambda} \in (0,+\infty]$ such that, for every $\lambda\in (0,\bar{\lambda})$, the problem \ref{Plambda} has at least two $L^p$-viscosity solutions, $u_{\lambda, 1}$ and $u_{\lambda , 2}\,$, satisfying
$$
u_{\lambda , 1} \xrightarrow[{\lambda\rightarrow 0^+}]{} u_0\;\; \textrm{in} \;E\,, \;\;\;\max_{\overline{\Omega}} u_{\lambda , 2} \xrightarrow[{\lambda\rightarrow 0^+}]{} +\infty\,,
$$
and, if $\,\bar{\lambda}<+\infty$, the problem $(P_{\bar{\lambda}})$ has at least one  $L^p$-viscosity solution. The latter is unique if $F(x,r,\vec{p},X)$ is convex in $(r,\vec{p},X)$.

4. If in addition \eqref{Hstrong} holds, the solutions $u_\lambda\,$ for $\lambda\leq 0$ are unique among $L^p$-viscosity solutions; whereas the solutions from 3. for $\lambda>0$ are ordered, $u_{\lambda, 1} \ll u_{\lambda , 2}$.

\end{teo}

This theorem proves the multiplicity conjectures from  \cite{arma2010}, \cite{note}.
We recall that Theorem \ref{th1.1,1.2,1.3} is new even when $F$ is linear operator in nondivergence form.

The supplementary hypotheses for the uniqueness results in the above theorem are unavoidable -- we recall that, in the universe of viscosity solutions, uniqueness is only available in the presence of a strong solution (see \cite{arma2010} and the references in that paper). \medskip

In the next two theorems, we show that it is possible to obtain a more precise description of the set $\Sigma$, provided we know the sign of $u_0$.  Such results for the divergence case $F=\Delta$ were already  proved in  theorems 1.4 and 1.5 in \cite{CJ}. Note that if $h$ has a sign, then $u_0$ has the same sign, by the maximum principle (see remark \ref{particularh}).

\begin{teo} \label{th1.5}
Suppose \ref{SC0}, \eqref{Hbeta}, \eqref{ExistUnic M bem definido}, \eqref{Hstrong} and \ref{H0} with $u_0 \leq 0$ and $cu_0\lneqq 0$. Then every nonpositive solution of \ref{Plambda} with $\lambda>0$ satisfies $u\ll u_0$. Furthermore, for every $\lambda >0$, the problem \ref{Plambda} has at least two nontrivial $L^p$-viscosity solutions $ u_{\lambda, 1} \ll u_{\lambda , 2}\,$, such that  $u_{\lambda_2 ,1}\ll u_{\lambda_1 ,1} \ll u_0$ if $\,0<\lambda_1<\lambda_2\,$, and
$$
u_{\lambda , 1} \xrightarrow[{\lambda\rightarrow 0^+}]{} u_0\;\; \textrm{in} \;E\,, \;\;\;\max_{\overline{\Omega}} u_{\lambda , 2} \xrightarrow[{\lambda\rightarrow 0^+}]{} +\infty\,.
$$
If  the operator $F(x,r,\vec{p},X)$ is convex in $(r,\vec{p},X)$ then $\max_{\overline{\Omega}}\, u_{\lambda,2}>0$ for all $\lambda>0$.
\end{teo}

\begin{figure}[H] \label{figurath15}
\centering
\includegraphics[scale=0.33]{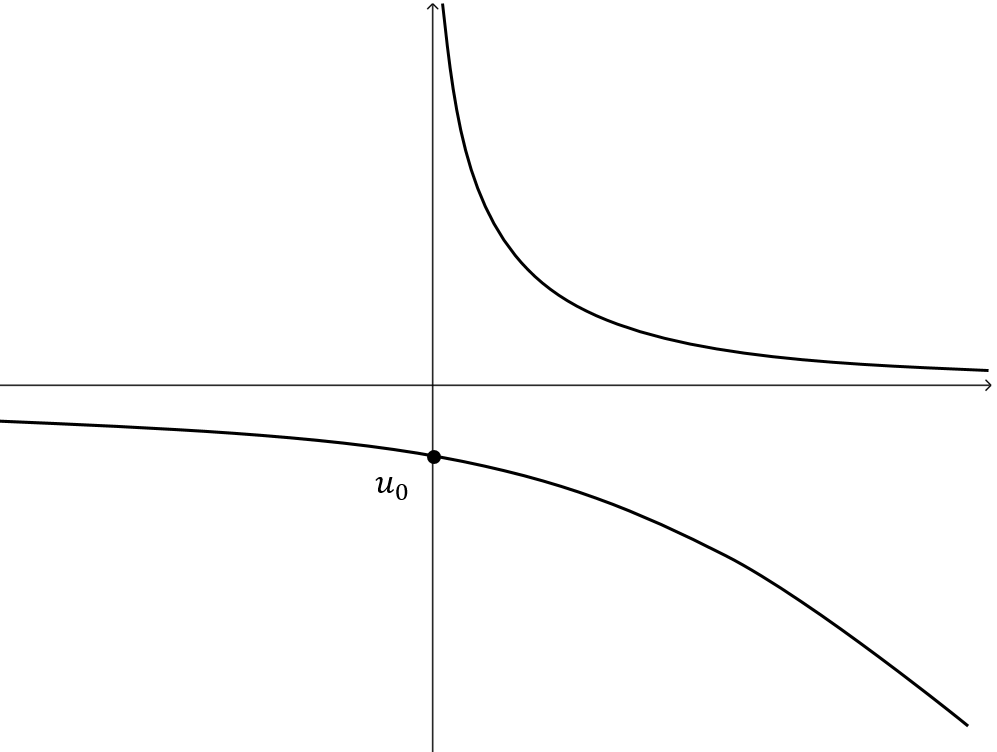}
\caption{Illustration of theorem \ref{th1.5}.}
\label{Rotulo1}
\end{figure}

In this figure we have put $\lambda$ on the horizontal axis. On the negative side of the vertical axis we have $u_{\lambda,1} (x_0)$ for any fixed $x_0\in\Omega$, which is a negative quantity for $\lambda>0$; whereas on the positive side of the vertical axis we find $\|u_{\lambda,2}\|_{L^\infty (\Omega)}$ (or $\max_{\overline{\Omega}} u_{\lambda,2}$ if $F$ is convex).

\begin{teo} \label{th1.4}
Suppose \ref{SC0}, \eqref{Hbeta}, \eqref{ExistUnic M bem definido}, \eqref{Hstrong} and \ref{H0} with $u_0 \geq 0$ and $cu_0\gneqq 0$. Then every nonnegative solution of \ref{Plambda} with $\lambda>0$ satisfies $u\gg u_0$. Moreover, there exists  $\bar{\lambda} \in (0,+\infty)$ such that
\begin{enumerate}[(i)]
\item for every $\lambda\in (0,\bar{\lambda})$, the problem \ref{Plambda} has at least two nontrivial $L^p$-viscosity solutions with $ u_{\lambda, 1} \ll u_{\lambda , 2}\,$, where $u_0 \ll u_{\lambda_1 ,1}\ll u_{\lambda_2 ,1} $ if $\,0<\lambda_1<\lambda_2\,$ and
$$
u_{\lambda , 1} \xrightarrow[{\lambda\rightarrow 0^+}]{} u_0\;\; \textrm{in} \;E\,, \;\;\;\max_{\overline{\Omega}} u_{\lambda , 2} \xrightarrow[{\lambda\rightarrow 0^+}]{} +\infty\,;
$$
\item the problem $(P_{\bar{\lambda}})$ has at least one $L^p$-viscosity solution $u_{\bar{\lambda}}$\,; this solution is unique if $F$ is convex in $(r,\vec{p},X)$;
\item for $\lambda > \bar{\lambda}$, the problem \ref{Plambda} has no nonnegative solution.
\end{enumerate}
\end{teo}

\begin{figure}[!htb]
\centering
\includegraphics[scale=0.33]{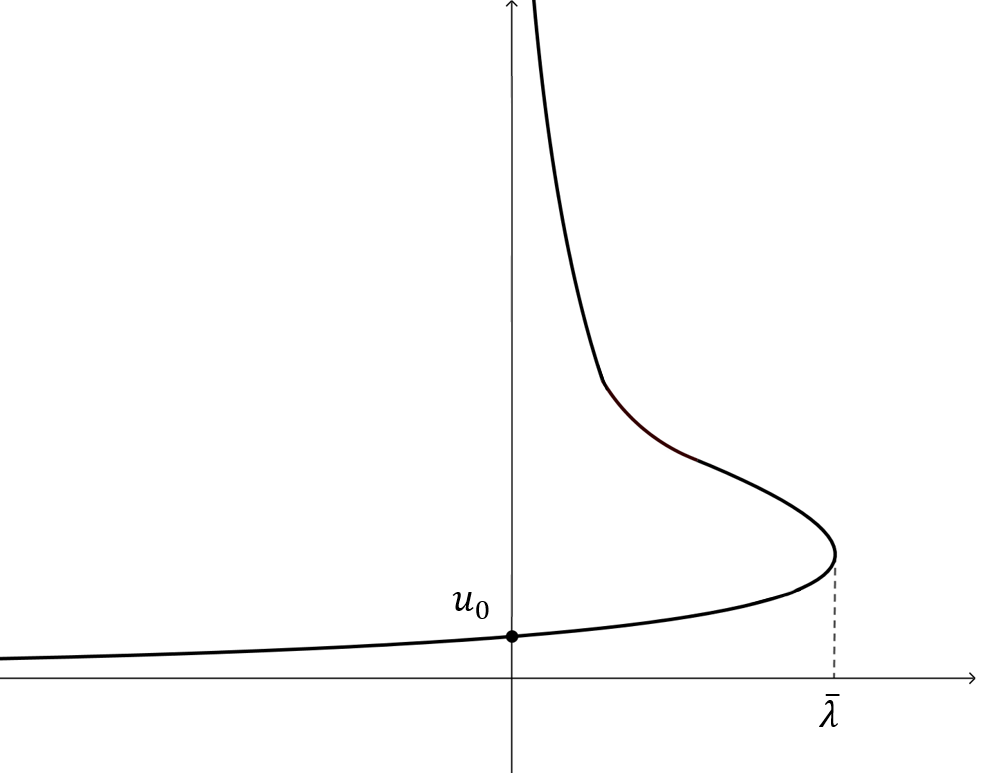}
\caption{Illustration of theorem \ref{th1.4}.}
\label{Rotulo2}
\end{figure}

\vspace{0.2cm}

Notice that, up to replacing $u$ for $-u$,  we are also taking into account, indirectly, the case $-\mu_2I\leq M(x)\leq -\mu_1 I\,$ for $\mu_1,\mu_2>0$. We only need to pay attention to the sign of $u_0$, which is reversed in this case.
\smallskip

In the end, we note that the hypothesis on the boundedness of the coefficients in \ref{SC0} is only due to the unavailability of any form of the Vázquez strong maximum principle for equations with unbounded coefficients, a rather interesting open problem in itself.


\section{Preliminaries}\label{Preliminaries}

Let $F(x,r,\vec{p},X):\Omega\times\real\times\rn\times\mathbb{S}^n\rightarrow\real$ be a measurable function satisfying \eqref{SC}. Notice that the condition on the zero order term in \eqref{SC} means that $F$ is proper, i.e. decreasing in $r$, while the hypothesis on the $X$ entry  implies that $F$ is a uniformly elliptic operator.
In \eqref{SC},
$$
\mathcal{M}^+(X):=\sup_{\lambda_P I\leq A\leq \Lambda_P I} \mathrm{tr} (AX)\,,\quad \mathcal{M}^-(X):=\inf_{\lambda_P I\leq A\leq \Lambda_P I} \mathrm{tr} (AX)
$$
are the Pucci's extremal operators with constants $0<\lambda_P\leq \Lambda_P$\footnote{We are denoting the ellipticity coefficients by $\lambda_P$ and $\Lambda_P$ instead of the usual $\lambda$ and $\Lambda$ in order to avoid any confusions with $\lambda$ in the problem \ref{Plambda}.}. See, for example, \cite{CafCab} for their properties.
Also, denote $\mathcal{L}^\pm [u]:=\mathcal{M}^\pm (D^2 u)\pm b(x)|Du|$, for $b\in L^p_+ (\Omega)$.

Now we recall the definition of $L^p$-viscosity solution.
\begin{defin}\label{def Lp-viscosity sol}
Let $f\in L^p_{\textrm{loc}}(\Omega)$.
We say that $u\in C(\Omega)$ is an $L^p$-viscosity subsolution $($respectively, supersolution$)$ of $F=f$ in $\Omega$ if whenever $\phi\in  W^{2,p}_{\mathrm{loc}}(\Omega)$, $\varepsilon>0$ and $\mathcal{O}\subset\Omega$ open are such that
\begin{align*}
F(x,u(x),D\phi(x),D^2\phi (x))-f(x)  \leq -\varepsilon\;\;
( F(x,u(x),D\phi(x),D^2\phi (x))-f(x)  \geq \varepsilon )
\end{align*}
for a.e. $x\in\mathcal{O}$, then $u-\phi$ cannot have a local maximum $($minimum$)$ in $\mathcal{O}$.
\end{defin}

We can think about $L^p$-viscosity solutions for any $p>\frac{n}{2}$, since this restriction makes all test functions $\phi\in W^{2,p}_{\mathrm{loc}}(\Omega)$ continuous and having a second order Taylor expansion \cite{CCKS}.
We are going to deal specially with the case $p>n$.


If $F$ and $f$ are continuous in $x$, we can use the more usual notion of  \textit{$C$-viscosity} sub and supersolutions -- see \cite{user}. Both definitions are equivalent when $F$ satisfies \ref{SC0} with $p\geq n$, by proposition 2.9 in \cite{CCKS}; we will be using them interchangeably, in this case, throughout the text.

On the other side, a \textit{strong} sub or subsolution belongs to $W^{2,p}_{\mathrm{loc}}(\Omega)$ and satisfies the inequality at almost every point. As we already mentioned, this notion is intrinsically connected to the $L^p$-viscosity concept; more precisely we have the following.
\begin{prop}\label{Lpiffstrong.quad}
Let $F$ satisfy \eqref{SC} and $f\in L^p(\Omega)$, $\mu\ge0$.
Then, $u\in W^{2,p}_{\mathrm{loc}}(\Omega)$ is a strong subsolution $($supersolution$)$ of $F+\mu|D u|^2=f$ in $\Omega$ if and only if it is an $L^p$-viscosity subsolution $($supersolution$)$ of this equation.
\end{prop}

See theorem 3.1 and proposition 9.1 in \cite{KSweakharnack} for a proof, even with more general conditions on $\mu$ and the exponents $p,q$.
When we refer to solutions of the Dirichlet problem, we will assume that strong solutions belong to $W^{2,p}(\Omega)$, unless otherwise specified.
Remember that a \textit{solution} is always both a sub and supersolution of the equation.
It is also well known that the pointwise maximum of subsolutions, or supremum over any set $($if this supremum is locally bounded$)$, is still a subsolution\footnote{See proposition 2 in \cite{KoikePerronRevisited} for a version for $L^p$-viscosity solutions, $p>n$, related to quadratic growth.}.

The next proposition follows from theorem 4 in \cite{arma2010} in the case $p=n$. For a version with more general exponents and coefficients we refer to proposition 9.4 in \cite{KSweakharnack}.

\begin{prop} {$($Stability$)$} \label{Lpquad}
Let $F$, $F_k$ be operators satisfying \eqref{SC}, $p\geq n$, $f, \, f_k\in L^p(\Omega)$, $u_k\in C(\Omega)$ an $L^p$-viscosity subsolution $($supersolution$)$ of
$$
F_k(x,u_k,Du_k,D^2u_k)+\langle M(x) Du_k,Du_k\rangle \geq(\leq) f_k(x) \;\;\textrm{in} \;\;\Omega\, , \;\textrm{ for all } k\in \n .
$$
Suppose $u_k\rightarrow u$ in $L_{\mathrm{loc}}^\infty  (\Omega)$ as $k\rightarrow \infty$ and, for each $B\subset\subset \Omega$ and $\varphi\in W^{2,p}(B)$, if we set
\begin{align*}
g_k(x):=F_k(x,u_k,D\varphi,D^2\varphi) \rangle- f_k(x) \, , \;
g(x):=F(x,u,D\varphi,D^2\varphi)-f(x)
\end{align*}
we have $\| (g_k-g)^+\|_{L^p(B)}$ $(\| (g_k-g)^-\|_{L^p(B)}) \rightarrow 0$ as $k\rightarrow \infty$. Then $u$ is an $L^p$-viscosity subsolution $($supersolution$)$ of\;
$
F(x,u,Du,D^2u)+\langle M(x)Du,Du\rangle  \geq(\leq) f(x) \,\textrm{ in }\, \Omega\, .
$
\end{prop}

The following result, which follows from lemma 2.3 in \cite{arma2010}, is a useful tool to deal with quadratic dependence in the gradient. Since we are going to use this result several times in the text, we present a proof in the appendix, for reader's convenience.

\begin{lem}{$($Exponential change$)$}\label{lemma2.3arma}
Let $p\geq n$ and $u\in W^{2,p}_{\mathrm{loc}}(\Omega)$. For $m>0$ set
$$
v=\frac{e^{mu}-1}{m} \, , \quad w=\frac{1-e^{-mu}}{m}.
$$
Then, a.e. in $\Omega$ we have $Dv=(1+mv)Du$, $Dw=(1-mw) Du$ and
\begin{align}
\mathcal{M}^\pm (D^2 u)+m\lambda_P |Du|^2 &\leq \frac{\mathcal{M}^\pm (D^2 v)}{1+mv} \leq \mathcal{M}^\pm (D^2 u)+m\Lambda_P |Du|^2 , \label{arma2010eq.uev}\\
\mathcal{M}^\pm (D^2 u)-m\Lambda_P |Du|^2 &\leq \frac{\mathcal{M}^\pm (D^2 w)}{1-mw} \leq \mathcal{M}^\pm (D^2 u)-m\lambda_P |Du|^2 .\label{arma2010eq.uew}
\end{align}
and clearly $\{u=0\}=\{v=0\}=\{w=0\}$ and $\{u>0\}=\{v>0\}=\{w>0\}$.

Moreover, the same inequalities hold in the $L^p$-viscosity sense if $u$ is merely continuous, that is, for example, if  $u\in C(\Omega)$ is an $L^p$-viscosity solution of
\begin{align}\label{arma2010eq.u}
\mathcal{M}^+(D^2 u)+b(x)|Du| + \mu |Du|^2 +c(x)u\geq f(x) \;\;\mathrm{in}\;\Omega
\end{align}
where $b,c,f\in L^p(\Omega)$, then $v=\frac{1}{m}(e^{mu} -1)$, for $m=\frac{\mu}{\lambda_P}$, is an $L^p$-viscosity solution of
\begin{align}\label{arma2010eq.v}
\mathcal{M}^+(D^2 v)+b(x)|Dv| +\frac{c(x)}{m} \mathrm{ln}(1+mv)(1+mv) -m f(x)v\geq f(x) \;\,\mathrm{in}\;\Omega
\end{align}
and analogously for the other inequalities.
\end{lem}

\vspace{0.1cm}
We recall some Alexandrov-Bakelman-Pucci type results with unbounded ingredients and quadratic growth, which will be referred to simply by ABP.

\begin{prop} \label{wABPquad}
Let $\Omega$ be bounded, $\mu\geq 0$ and $f\in L^p (\Omega)$ for $p>n$.
Then, there exist $\delta=\delta (n,p,\lambda_P, \Lambda_P, \mathrm{diam}(\Omega), \|b\|_{L^p(\Omega)})>0$ such that if
\begin{align*}
\quad\quad\quad\quad \mu \|f^-\|_{L^p(\Omega)} \, (\mathrm{diam}(\Omega))^{\frac{n}{p}} \leq \delta  \quad\quad\quad (\mu \|f^+\|_{L^p(\Omega)}  \, (\mathrm{diam}(\Omega))^{\frac{n}{p}} \leq \delta )
\end{align*}
then every $u\in C(\overline{\Omega})$ which is an $L^p$-viscosity subsolution $($supersolution$)$ of
\begin{align*}
\mathcal{L}^+[u]+\mu |D u|^2 \geq f(x)\;\; \mathrm{in}\;\; \Omega \qquad
\left( \mathcal{L}^-[u]-\mu |D u|^2 \leq f(x)\;\; \mathrm{in}\;\; \Omega  \,\right)
\end{align*}
satisfies, for a constant $C=C (n,p,\lambda_P,\Lambda_P, \|b\|_{L^n(\Omega)},\mathrm{diam}(\Omega))$, the estimate
\begin{align*}
\max_{\overline{\Omega}} u \leq \max_{\partial \Omega} u +C \,  \|f^-\|_{L^p(\Omega)} \quad \left( \min_{\overline{\Omega}} u \geq \min_{\partial \Omega} u -C  \|f^+\|_{L^p(\Omega)} \right).
\end{align*}
\end{prop}

We refer to theorem 2.6 of \cite{KSweakharnack} (see also \cite{KSexist2009}, \cite{KSmpite2007}), for more general forms of ABP. 

A consequence of ABP in its quadratic form is the  comparison principle for our equations, concerning $L^p$-viscosity solutions and coercive operators. We give a short proof for reader's convenience.
We make the convention, along the text, that $\alpha$ and $\beta$ will always denote a pair of sub and supersolutions, in a sense to be specified.

\begin{lem} \label{lema subsol strong always less supersol}
Assume $F$ satisfies \eqref{SC}, $f\in L^p(\Omega)$ and $M\in L^\infty (\Omega)$. Suppose that $u$ is an $L^p$-viscosity supersolution of
\begin{align}\label{P com f e sem c}
\left\{
\begin{array}{rclcc}
-F[u]&=& \langle M(x) D u,D u\rangle +f(x) & \mbox{in} & \Omega\\
u &=& 0 & \mbox{on} & \partial\Omega
\end{array}
\right.
\end{align}
Then, for any $\alpha\in C(\overline{\Omega})\cap W^{2,p}_{\mathrm{loc}}(\Omega)$ strong subsolution of \eqref{P com f e sem c}, we have $ \alpha\leq u$ in $\Omega$.

Analogously, if $u$ is an $L^p$-viscosity subsolution of \eqref{P com f e sem c}, for any $\beta\in C(\overline{\Omega})\cap W^{2,p}_{\mathrm{loc}}(\Omega)$ strong supersolution of \eqref{P com f e sem c}, we have $ u\leq \beta$ in $\Omega$.
\end{lem}

\begin{proof}
Set $v:=u-\alpha$ in $\Omega$. By contradiction, assume
$
\min_{\overline{\Omega}} v = v(x_0)<0.
$
As $v\geq 0$ on $\partial\Omega$, thus $x_0\in\Omega$. Set $\widetilde{\Omega}:=\{v<0 \textrm{ in }\Omega \}$, which is an open nonempty set since $x_0\in\widetilde{\Omega}$. Let $\varphi\in W^{2,p}_{\mathrm{loc}}(\widetilde{\Omega})$ and $\tilde{x}\in\widetilde{\Omega}$ be such that $v-\varphi$ has a minimum at $\tilde{x}$. But then $u-(\alpha+\varphi)$ has a minimum at $\tilde{x}$, and by $\alpha+\varphi\in W^{2,p}_{\mathrm{loc}}(\widetilde{\Omega})$ together with the definition of $u$ being an $L^p$-viscosity supersolution, we know that
for every $\varepsilon >0$, there exists $r>0$ such that, for a.e. $x\in B_r(\tilde{x})\cap\widetilde{\Omega}$, we have
$$
- F(x,u,D\alpha+D\varphi,D^2\alpha+D^2\varphi)-\langle M(x)D (\alpha+\varphi),D (\alpha+\varphi)\rangle -f(x)\geq -\varepsilon
$$
and $- F(x,\alpha,D\alpha,D^2\alpha)-\langle M(x)D \alpha,D \alpha\rangle -f(x)\leq 0$ from the definition of $\alpha$.
By \eqref{SC},
\begin{align*}
\varepsilon &\geq F(x,u,D\alpha+D\varphi,D^2\alpha+D^2\varphi)- F(x,\alpha,D\alpha,D^2\alpha)\\
&\hspace{0.5cm}+\langle M(x)D (\alpha+\varphi),D (\alpha+\varphi)\rangle-\langle M(x)D \alpha,D \alpha\rangle \\
&\geq \mathcal{M}^-(D^2\varphi)-b(x)|D\varphi|-\mu\, (\,|D\alpha+D\varphi|+|D\varphi|\,)\,|D\varphi| -d(x)\,\omega((u-\alpha)^+) \\
&\geq \mathcal{M}^-(D^2\varphi)-b(x)|D\varphi|-\mu |D\varphi|^2-2\mu |D\alpha|\,|D\varphi|
\end{align*}
which means that, for $\widetilde{b}=b+2\mu |D\alpha|\in L^p_+(\Omega)$, $v$ is an $L^p$-viscosity supersolution of
\begin{align*}
\left\{
\begin{array}{rclcl}
\mathcal{M}^-(D^2 v)-\widetilde{b}(x)|Dv|-\mu |Dv|^2 &\leq & 0 &\mbox{in}& \;\widetilde{\Omega}\\
v &\geq & 0 &\mbox{on}& \partial\widetilde{\Omega}\subset\partial\Omega\cup\{ v=0\}.
\end{array}
\right.
\end{align*}
Then, by ABP with $f=0$, we have $v\geq 0$ in $\widetilde{\Omega}$, which contradicts the definition of $\widetilde{\Omega}$.
\end{proof}

\begin{rmk}
The same result hods if $\,\alpha =\max_{1\leq i\leq m}\alpha_i$ and $\,\beta = \min_{1\leq j\leq l}\beta_j\,$, where $\alpha_i$ and $\beta_j$ are continuous strong sub and supersolutions of \eqref{P com f e sem c} respectively. Indeed, in the proof above, we only need to note that \,$\min_{\overline{\Omega}}\,(u-\alpha)=(u-\alpha_i) (x_0)\,$ for some $i\in \{1, \ldots , m\}$, and consider $v:=u-\alpha_i\,$.
\end{rmk}

\begin{rmk} \label{obs.c.leq0}
By lemma \ref{lema subsol strong always less supersol} we obtain that for any pair of $L^p$-viscosity sub and supersolutions $\alpha,\beta$ of $(P_0)$, we have $\alpha\leq u\leq \beta$ in $\Omega$, for any strong solution $u$ of $(P_0)$.
\end{rmk}

\vspace{0.05 cm}

The Local Maximum Principle (LMP) is well known, see for example \cite{GT}, \cite{CafCab}, \cite{Koike}, and \cite{KSlmp}. Its  boundary version (BLMP) is given in theorem 1.3 in \cite{B2016}, without zero order term.   We denote with $B_R^+=\{x\in\rn\::\:|x|<R,\;x_n>0\}$ a half ball with a flat portion of the boundary included in $\{x_n=0\}$.

\begin{prop} \label{BLMP}
Let $u$ be a locally bounded $L^p$-viscosity subsolution of
\begin{align*}
\left\{
\begin{array}{rclcl}
\mathcal{L}^+ (D^2 u)+\nu (x) u &\geq & -f(x) &\mbox{in} & B_2^+ \\
u &\leq & 0 &\mbox{on} & B_{2}\cap \{x_n=0\}
\end{array}
\right.
\end{align*}
with $f\in L^p (B_2^+)$, $b\in L^q_+(\Omega)$, $\nu\in L^{p_1} (B_2^+)\cap L^p (B_2^+)$, for some $p_1>n$ and $q\geq p\geq n$, $q>n$.
Then, for each $r>0$,
\begin{align*}
 \sup_{B_1^+} u \leq C \left( \left(  \int_{B_{3/2}^+} (u^+)^r \right)^{1/r} + \|f^+\|_{L^p(B_2^+)} \right),
\end{align*}
where $C$ depends only on $n,\,p, p_1,\,\lambda, \,\Lambda,\, r,\, \|b\|_{L^q(B_2^+)}$ and $\,\|\nu\|_{L^{p_1}(B_2^+)}$.
\end{prop}

The proof, as in \cite{B2016}, is a consequence of LMP. The only difference comes from the need to put the zero order term on the right hand side, applying theorem 3.1(a) of \cite{KSlmp}, followed by a Moser type approach as in \cite{HanLin}; a detailed proof is given in \cite{tese}.


We recall two boundary versions of the quantitative strong maximum principle and the weak Harnack inequality, which follow by theorems 1.1. and 1.2 in \cite{B2016}.

\begin{teo}\label{BQSMP}
Suppose $d,f\in L^p(B_2^+)$, $p>n$. Assume that $u$ is an $L^p$-viscosity supersolution of $\mathcal{L}^-[u]-du\leq f$, $u\geq 0$ in $B_2^+$. Then there exist constants $\varepsilon, c, C>0$ depending on $n,\lambda, \Lambda, p$, $\|b\|_{L^p (B_2^+)}$ and $\|d\|_{L^p (B_2^+)}$ such that
$$
\inf_{B_1^+} \frac{u}{x_n}\geq c \left( \int_{B_{3/2}^+} (f^-)^\varepsilon \right)^{1/\varepsilon} -C \|f^+\|_{L^p (B_2^+)}
$$
\end{teo}

\begin{teo}\label{BWHI}
Suppose $d,f\in L^p(B_2^+)$, $p>n$. Assume that $u$ is an $L^p$-viscosity supersolution of $\mathcal{L}^-[u]-du\leq f$, $u\geq 0$ in $B_2^+$. Then there exist constants $\varepsilon, c, C>0$ depending on $n,\lambda, \Lambda, p$, $\|b\|_{L^p (B_2^+)}$ and $\|d\|_{L^p (B_2^+)}$ such that
$$
\inf_{B_1^+} \frac{u}{x_n}\geq c \left( \int_{B_{3/2}^+} \left(\frac{u}{x_n} \right)^\varepsilon \right)^{1/\varepsilon} -C \|f^+\|_{L^p (B_2^+)}
$$
\end{teo}

Theorem \ref{BWHI} implies, in particular, the strong maximum principle when $f=0$, i.e. for $\Omega\in C^{1,1}$ and $u$ an $L^p$-viscosity solution of $\mathcal{L}^-[u]-du\leq 0$, $u\geq 0$ in $\Omega$, where $d\in L^p(\Omega)$, we have either $u\equiv 0$ in $\Omega$ or $u>0$ in $\Omega$ and if $u(x_0)=0$ at $x_0\in\partial\Omega$, then $\partial_\nu u (x_0)>0$. We will refer to these consequences, along the text, simply by SMP and Hopf.

In \cite{B2016}, theorems \ref{BQSMP} and \ref{BWHI} are proved for $d\equiv 0$, but exactly the same proofs there work for any $d\ge0$. Moreover, since the function $u$ has a sign, $d^-u\geq 0$ and they are also valid for nonproper operators.

The next two propositions are basic results in the study of first eigenvalues of weighted operators. For a proof with unbounded coefficients, see \cite{regularidade}.

\begin{prop}\label{exist eig for L-c geq 0}
Let $c\in L^p(\Omega)$, $c\gneqq 0$ for $p>n$. Then $\mathcal{L}^-$ has a positive weighted eigenvalue $\lambda_1>0$ corresponding to a positive eigenfunction $\varphi_1\in W^{2,p}(\Omega)$ such that
\begin{align} \label{eq exist eigen L-c}
\left\{
\begin{array}{rclcc}
(\mathcal{L}^-+\lambda_1c)\, [\varphi_1]&=&0 &\mbox{in} & \Omega\\
\varphi_1 &>& 0 &\mbox{in} & \Omega \\
\varphi_1 &=&0 &\mbox{on} &\partial\Omega.
\end{array}
\right.
\end{align}
\end{prop}

\begin{prop}\label{th4.1 QB}
Let $u,v\in C(\overline{\Omega})$ be $L^p$-viscosity solutions of
\begin{align}\label{def phi1- of L+(c)}
\begin{cases}
\mathcal{L}^-[u]+c(x)u\leq 0 \;&\mbox{in}\;\;\Omega\\
\hspace{2.05cm} u>0\;&\mbox{in}\;\;\Omega
\end{cases} \qquad\mbox{and}\qquad
\begin{cases}
\mathcal{L}^-[v]+c(x)v\geq 0 \;&\mbox{in}\;\;\Omega\\
\hspace{2cm} v\leq0\;&\mbox{on} \;\;\partial\Omega\\
\hspace{1.3cm} v(x_0)>0 &x_0 \in\Omega,
\end{cases}
\end{align}
with $c\in L^p(\Omega)$ for $p>n$. Suppose one, $u$ or $v$, is a strong solution, i.e. belongs to the space $W^{2,p}_{\mathrm{loc}}(\Omega)$.
Then $u=tv$ for some $t>0$.
\end{prop}

Proposition \ref{th4.1 QB} in \cite{regularidade} was stated in terms of $L^n$-viscosity solutions. However, for equations satisfying \eqref{SC}, the notions of $L^p$ and $L^n$ viscosity solutions are equivalent, by theorem 2.3(iii) in \cite{CKSS} (by using in its proof the ABP with unbounded coefficients from proposition 2.8 in \cite{KSmpite2007}); the addition of the term $c(x)u$ does not change the proof.


\section{Existence results through fixed point theorems} \label{secao fixed point problem}

Consider the problem \ref{Plambda} without $\lambda$ dependence, i.e.
\begin{align}\label{P} \tag{P}
\left\{
\begin{array}{rclcc}
-F(x,u,Du,D^2u)&=&c(x)u+\langle M(x)D u, D u \rangle +h(x) &\mbox{in} & \Omega\\
u &=& 0 &\mbox{on} & \partial\Omega
\end{array}
\right.
\end{align}
under \eqref{SC}.
In this section all results hold for functions $b$, $c$ and $h$ in $ L^p(\Omega)$.
About the matrix $M$, it  only need to assume that $M\in L^\infty (\Omega)$. We set $\mu = \|M\|_{L^\infty(\Omega)}$. As for $c$, no sign condition is needed in this section.

We define, under hypothesis \eqref{ExistUnic M bem definido}, the operator $\mathcal{T}: \, E\rightarrow E$ that takes a function $u\in E=C^1(\overline{\Omega})$ into $U=\mathcal{T}u$\,, the unique $L^p$-viscosity solution of the problem
\begin{align}\label{T_u} \tag{$\mathcal{T}_u$}
\left\{
\begin{array}{rclcc}
-F(x,U,DU,D^2U)&=& c(x)u+\langle M(x)D u, D u \rangle +h(x) &\mbox{in} & \Omega\\
U &=& 0 &\mbox{on} & \partial\Omega
\end{array}
\right.
\end{align}

\begin{claim} \label{T is completely continuous}
The operator $\mathcal{T}$ is completely continuous.
\end{claim}

\begin{proof}
 Let $u_k\in E$, $u_k\rightarrow u$ in $E$. Then $\|u_k\|_{E}\leq C_0$, for all $k\in\n$. Set
$f_k(x):=c(x)u_k+\langle M(x)D u_k, D u_k \rangle +h(x)$
so
$$
\|f_k\|_{L^p(\Omega)} \leq \|h\|_{L^p(\Omega)} + C_0 \,\|c\|_{L^p(\Omega)} + \mu\, C_0^2 \leq C.
$$
We use ABP on the  sequence $U_k=\mathcal{T}{u_k}$, from where we get $\|U_k\|_{L^\infty}\leq C \,\|f_k\|_{L^p(\Omega)}\leq C$. Thus the $C^{1,\alpha}$ global estimate from \cite{regularidade} gives us $\|U_k\|_{C^{1,\alpha}(\overline{\Omega})}\leq C .$
By the compact inclusion of $C^{1,\alpha} (\overline{\Omega})$ into $E$, there exists $U\in E$ and a subsequence such that $U_k \rightarrow U$ in $E$.
This already shows that $\mathcal{T}$ takes bounded sets into precompact ones.

We need to see that $U=\mathcal{T}u$. This easily follows from  the stability proposition \ref{Lpquad}, by defining, for each $\varphi\in W^{2,p}_{\mathrm{loc}}(\Omega)$,
\begin{align*}
g_k(x):= -F(x,U_k,D\varphi,D^2\varphi) - f_k(x)\,  , \;\;\;
g(x):=-F(x,U,D\varphi,D^2\varphi) - f(x)
\end{align*}
a.e. $x\in\Omega$, where $f$ is the same as $f_k$ with  $u_k$ replaced by $u$. Indeed, by \eqref{SC},
\begin{align*}
\|g_k&-g\|_{L^p(\Omega)} \leq \|d\|_{L^p(\Omega)} \,\omega(\|U_k-U\|_{L^\infty(\Omega)})  +\|c\|_{L^p(\Omega)} \|u_k-u\|_{L^\infty(\Omega)} \\
&+ \mu\, (\|D u_k\|_{L^\infty(\Omega)}+\|D u\|_{L^\infty(\Omega)})\,\|D u_k-D u\|_{L^\infty(\Omega)} |\Omega|^{1/p}\\
& \leq \|d\|_{L^p(\Omega)} \,\omega(\|U_k-U\|_{L^\infty(\Omega)}) +\{\,\|c\|_{L^p(\Omega)} +2 C_0 \,\mu\, |\Omega|^{1/p} \} \,\|u_k-u\|_{E}
\xrightarrow[{k\rightarrow +\infty}]{} 0\, ,
\end{align*}
since $\omega$ is increasing, continuous in $0$, with $\omega(0)=0$. Since the problem \eqref{T_u} has a unique solution, $U=\mathcal{T}u$.
On the other hand, since this argument can be made for any subsequence of the original $(U_k)_{k\in\n}$,  the whole sequence converges to $U=\mathcal{T}u$.

\qedhere{\,\textit{Claim \ref{T is completely continuous}.}}
\end{proof}

The next existence statement is a typical result about existence between sub and supersolutions, and is a version of theorem 2.1 of \cite{CJ} for fully nonlinear equations. We start with a definition.
\begin{defin} \label{def2.1}
An $L^p$-viscosity subsolution $\alpha\in E$ $\mathrm{(}$respectively, supersolution $\beta)$ of $(\mathrm{P})$ is said to be strict if every $L^p$-viscosity supersolution $($subsolution$)$ $u\in E$ of $(\mathrm{P})$ such that $\alpha\leq u$ $(u \leq\beta)$ in $\Omega$, also satisfies $\alpha\ll u$ $(u\ll \beta)$ in $\Omega$.
\end{defin}

Set $B_r=B_r^{E}(0):= \{u\in E; \|u\|_{E}<r\}$, for any $r>0$.
\begin{teo}\label{th2.1}
Let $\Omega \subset \rn$ be a bounded domain with $\partial\Omega\in C^{1,1}$. Suppose \eqref{SC}, \eqref{Hbeta} and \eqref{ExistUnic M bem definido}.
Let $\alpha =\max_{1\leq i \leq \kappa} \,\alpha_i \,,\, \beta = \min_{1\leq j \leq \iota} \,\beta_j \,$, where $\alpha_i \, , \,\beta_j \in W^{2,p} (\Omega)$ are strong sub and supersolutions of \eqref{P} respectively, with $\alpha \leq \beta$ in $\Omega$.
Then \eqref{P} has an $L^p$-viscosity solution satisfying $\alpha \leq u \leq \beta$ in $\Omega$.
Furthermore,

\begin{enumerate}[(i)]
\item If $\alpha$ and $\beta$ are strict in the sense of definition \ref{def2.1}, then for large $R>0$ we have
$$
deg(I-\mathcal{T}, \mathcal{S},0)=1
$$
where $\mathcal{S}=\mathcal{O}\cap B_R$, for
$
\mathcal{O}=\mathcal{O}_{\alpha,\beta}:= \{u\in C_0^1(\overline{\Omega}); \; \alpha\ll u \ll \beta \;in \; \Omega\}.
$
\item If \eqref{Hstrong} holds, there exists a minimal solution $\underline{u}$ and a maximal solution $\overline{u}$  of \eqref{P} in the sense that every $($strong$)$ solution $u$ of \eqref{P} in the order interval $[\alpha , \beta]$ $($i.e. such that $\alpha (x)\leq u(x)\leq \beta (x)$ for all $ x\in\Omega)$ satisfies
$$
\alpha\,\leq \,\underline{u} \,\leq u \leq \,\overline{u} \,\leq \beta \;\;\; \mathrm{in} \;\;\; \Omega.
$$
\end{enumerate}
\end{teo}

Notice that, under the assumptions of theorem \ref{th2.1} and by the  global $C^{1,\alpha}$-estimate \cite{regularidade}, every solution $u$ of \eqref{P} satisfying $\alpha\leq u\leq \beta$ in $\Omega$ is such that
\begin{align}\label{u C1,alpha leq C}
\|u\|_{C^{1,\alpha}(\overline{\Omega})}\leq C,
\end{align}
where $C$ depends on $r_0,n,p,\lambda_P,\Lambda_P,\mu,\| b \|_{L^p (\Omega)}, \| c \|_{L^p (\Omega)}$, $\| h \|_{L^p (\Omega)}$, $\omega (1)\| d \|_{L^p (\Omega)}$, $\mathrm{diam} (\Omega)$ and, of course, on the $L^\infty$ uniform bounds on $u$ given by $\|\alpha\|_\infty $ and $ \|\beta\|_\infty$.

\begin{proof}
Consider any $R\geq \max \{C,\|\alpha\|_{E},\|\beta\|_{E} \}+1$, with $C$ from \eqref{u C1,alpha leq C}.

\textit{Part 1. Existence of a solution in the order interval $[\alpha,\beta]$.}

First of all, we construct a modified problem, similar to but a little bit simpler than the one given in \cite{CJ}. In order to avoid technicalities, consider $\kappa=\iota=1$; later we indicate the corresponding changes. We set
$f(x,r,\vec{p})= h(x)+c(x)\,r+\langle M(x)\,\vec{p},\vec{p} \rangle$; $\overline{f}(x,r,\vec{p})=  h(x)+c(x)\,r+\langle \overline{M}(x)\,\vec{p},\vec{p} \rangle$ for
\begin{align*}
\overline{M}(x)=\overline{M}(x,\vec{p}):=
\begin{cases}
\;M(x) \, , \quad\quad\;\;\; \mathrm{if}\;\; |\vec{p}|<R \\
\;M(x) \frac{R^2}{|\vec{p}|^2}\, , \quad\;\; \mathrm{if}\;\; |\vec{p}|\geq R,
\end{cases}
\end{align*}
and also
\begin{align*}
\widetilde{f}(x,r,\vec{p}):=
\begin{cases}
\;\;\overline{f}(x,\alpha (x),D\alpha (x))\, , & \mathrm{if}\;\; r<\alpha (x) \\
\;\;\;\overline{f}(x,r,\vec{p}) \, ,& \mathrm{if}\;\;\; \alpha\leq r \leq \beta (x) \\
\;\;\overline{f}(x,\beta (x),D\beta (x))\, , & \mathrm{if}\;\; r>\beta (x).
\end{cases}
\end{align*}
Consider the problems
\begin{align}
\tag{P}
-F[u] =f(x,u,D u) \;\;\mbox{ in }  \Omega, \\
\label{Pbar} \tag{$\overline{\mathrm{P}}$}
-F[u]=\overline{f}(x,u,D u) \;\;\mbox{ in } \Omega,\\
\label{Ptilde} \tag{$\widetilde{\mathrm{P}}$}
-F[u]=\widetilde{f}(x,u,D u) \;\;\mbox{ in }  \Omega,\\
\mbox{and } u=0\;\;\mbox{ on }  \partial\Omega.\nonumber
\end{align}

Notice that, by the construction of $\overline{f}, \, \widetilde{f}$ and $|D\alpha|, \, |D\beta|<R$, we have that
\begin{align}\label{ftilde=f}
\widetilde{f}(x,\alpha (x),D\alpha (x))=\overline{f}(x,\alpha (x),D\alpha (x))=f(x,\alpha (x),D\alpha (x)) \;\;\; \mathrm{a.e.\;in\;} \Omega
\end{align}
and the same for $\beta$. So, $\alpha,\beta$ are also a pair of strong sub and supersolutions for ($\overline{\mathrm{P}}$) and ($\widetilde{\mathrm{P}}$).
Observe also that, since $\|\overline{M}\|_{L^\infty(\Omega)}\leq\mu\,$ (and the operator $F(x,r,p,X)+\langle \overline{M}(x)p,p\rangle$ still satisfies the conditions in \cite{regularidade}), by \eqref{u C1,alpha leq C} and the definition of $R$, every solution $u$  of ($\overline{\mathrm{P}}$) with $\alpha\leq u \leq \beta$ in $\Omega$ is $C^{1,\alpha}$ up to the boundary and satisfies
\begin{align} \label{u<R}
\|u\|_{E}<R\,,
\end{align}
so $\overline{M}(x)=\overline{M}(x,D u)=M(x)$, and $u$ is a solution of the original problem (P).
\begin{claim} \label{claim1}
Every $L^p$-viscosity solution $u$ of $(\widetilde{\mathrm{P}})$ satisfies $\alpha\leq u \leq\beta$ in $\Omega$, hence is a solution of $(\overline{\mathrm{P}})$, and, by the above, a solution of $(\mathrm{P})$.
\end{claim}
\begin{proof}
Let $u$ be an $L^p$-viscosity solution of ($\widetilde{\mathrm{P}}$). As in the proof of of lemma \ref{lema subsol strong always less supersol}, in order to obtain a contradiction, suppose that $v:=u-\alpha$ is such that
$
\min_{\overline{\Omega}} v = v(x_0)<0\,.
$
As $v\geq 0$ on $\partial\Omega$, it follows that $x_0\in\Omega$. Consider $\widetilde{\Omega}:=\{ v<0 \}\neq\emptyset$.

We claim that $v$ is an $L^p$-viscosity supersolution of
\begin{align} \label{vsupersol}
\mathcal{M}^-(D^2 v)-b(x)|Dv| \leq 0 \;\textrm{ in }  \widetilde{\Omega}.
\end{align}
Indeed, let $\varphi\in W^{2,p}_{\mathrm{loc}}(\widetilde{\Omega})$ and $\widehat{x}\in\widetilde{\Omega}$ such that $v-\varphi$ has a local maximum at $\widehat{x}$ in $\widetilde{\Omega}$. Then, $\alpha+\varphi\in W^{2,p}_{\mathrm{loc}}(\widetilde{\Omega})$ and $u-(\alpha+\varphi)$ has also a local maximum at $\widehat{x}$.
By definition of $u$ as an $L^p$-viscosity supersolution of ($\widetilde{\mathrm{P}}$), for every $\varepsilon>0$, there exists an $r>0$ such that
\begin{align}\label{eq acima}
-F(x,u,D(\alpha+\varphi,D^2(\alpha+\varphi)))-\overline{f}(x,\alpha,D \alpha) \geq -\varepsilon \;\textrm{ a.e. in } B_r(\widehat{x})\cap\widetilde{\Omega}\,,
\end{align}
since $u<\alpha$ in $\widetilde{\Omega}$. Further, from \eqref{ftilde=f},
$-F(x,\alpha,D\alpha,D^2\alpha)-\overline{f}(x,u,D \alpha) \leq 0$ a.e. in $\widetilde{\Omega}$.
Subtracting this from \eqref{eq acima} and using \eqref{SC}, we obtain
$\mathcal{M}^- (D^2 \varphi)-b(x)|D\varphi|\leq \varepsilon$ a.e. in $B_r(\widehat{x})\cap\widetilde{\Omega}$,
since $d(x)\,\omega(v^+) \equiv0$ in $\widetilde{\Omega}$. Hence, \eqref{vsupersol} is proved.

Now, \eqref{vsupersol} and ABP imply that $v\geq 0$ in $\widetilde{\Omega}$ , which contradicts the definition of $\widetilde{\Omega}$. So, $\min_{\overline{\Omega}}v>0$ i.e. $u\geq \alpha$ in $\Omega$. Analogously, $u\leq\beta$ in $\Omega$.
\qedhere{\textit{Claim \ref{claim1}.}}
\end{proof}

Next, we move to build a solution for the operator $\widetilde{\mathcal{T}}: \, E\rightarrow E$ that takes a function $u\in E$ into $U=\widetilde{\mathcal{T}}_u$\,, the unique $L^p$-viscosity solution of the problem
\begin{align}\label{Ttilde_u} \tag{$\widetilde{\mathcal{T}}_u$}
\left\{
\begin{array}{rclcc}
-F[U]&=&\widetilde{f}(x,u,D u) &\;\;\mbox{in} & \Omega\\
U &=& 0 &\;\;\mbox{on} & \partial\Omega.
\end{array}
\right.
\end{align}

Solutions of ($\widetilde{\mathrm{P}}$) are fixed points of $\widetilde{\mathcal{T}}$, and belong to the order interval $[\alpha,\beta]$, by claim \ref{claim1}. Moreover,
\begin{align} \label{cotaTtilde}
\|\widetilde{\mathcal{T}} u\|_{E} <R_0\, , \;\; \textrm{ for all }\, u \in E,
\end{align}
for an appropriate $R_0>R$. In fact, by observing that
\begin{align*}
|\langle\,\overline{M}(x)\,\vec{p},\vec{p}\rangle |=
\begin{cases}
\,|\langle M(x)\,\vec{p},\vec{p} \rangle | \, , \quad\quad\;\;\; \mathrm{if}\;\; |\vec{p}|<R \\\,\frac{R^2}{|\vec{p}|^2} |\langle M(x)\,\vec{p},\vec{p} \rangle | \, , \quad\;\; \mathrm{if}\;\; |\vec{p}|\geq R
\end{cases}
\leq \;\,\mu\, R^2 \, , \;\textrm{ for all } \vec{p}\in \rn,
\end{align*}
then
$
|\widetilde{f} (x,r,\vec{p})|\leq |h(x)|+\mu\, R^2+|c(x)| \max \{\|\alpha\|_\infty , \|\beta\|_\infty \}
$
and if $\gamma (x):=\widetilde{f} (x,u(x),D u (x))$,
$$
\|\gamma\|_{L^p(\Omega)} \leq \|h\|_{L^p(\Omega)} + \mu\, R^2 +\|c\|_{L^p(\Omega)} \max \{\|\alpha\|_\infty , \|\beta\|_\infty \},
$$
for every $u\in E$. Thus, the $C^{1,\alpha}$-estimates from \cite{regularidade} and ABP, applied to the problem \eqref{Ttilde_u}, give us that $U=\widetilde{\mathcal{T}} u \in C^{1,\alpha}(\overline{\Omega})$ and $\|U\|_{C^{1,\alpha} (\overline{\Omega})} \leq C <R_0$
and \eqref{cotaTtilde} follows.

Notice that $\widetilde{\mathcal{T}}$ takes bounded sets in $E$ into precompact ones, by the above and the compact inclusion $C^{1,\alpha}\subset E$. Also, if $u_k\rightarrow u$ in $E$, then $U_k=\widetilde{\mathcal{T}} u_k\rightarrow U$ in $E$ up to a subsequence; thus we can conclude that $\widetilde{\mathcal{T}}$ is completely continuous if we show that $U=\widetilde{\mathcal{T}} u$. To prove the latter, similarly to the argument with $\mathcal{T}$ we set, for $\varphi\in W^{2,p}_{\mathrm{loc}}(\Omega)$,
\begin{align*}
g_k(x):=& -F(x,U_k,D\varphi, D^2 \varphi)-c(x)\,\widetilde{u}_k-\langle\, \overline{M}(x,Du_k)\,D u_k,D u_k \rangle -h(x) \\
g(x):=& -F(x,U,D\varphi, D^2 \varphi)-c(x)\,\widetilde{u}-\langle\, \overline{M}(x,Du)\,D u,D u  \rangle -h(x),
\end{align*}
where, for each fixed function $v$, the tilde function $\widetilde{v}$ is the following truncation
\begin{align*}
\widetilde{v}(x):=
\left\{
\begin{array}{rcl}
\alpha (x) & \mbox{if} & v(x)<\alpha (x)\\
v (x) & \mbox{if} &\alpha (x)\leq v(x) \leq \beta (x)\\
\beta (x) & \mbox{if} &v(x)>\beta (x).
\end{array}
\right.
\end{align*}
Observe that, as it is elementary to check,
$\|\widetilde{u}_k - \widetilde{u}\|_{L^\infty (\Omega)}\leq \|u_k - u\|_{L^\infty (\Omega)} \rightarrow 0 $ as $k\rightarrow\infty$.
By the estimates for $\mathcal{T}$ and using that the function $\vec{p}\mapsto \langle \overline{M}(x,\vec{p})\vec{p},\vec{p}\rangle$ is continuous in $\vec{p}$, we get
$\|g_k-g\|_{L^p(\Omega)}  \rightarrow0$ as $k\rightarrow\infty$ and
so $U=\widetilde{\mathcal{T}}u$.

From the complete continuity of $\mathcal{T}$ and \eqref{cotaTtilde}, the degree $\mathrm{deg} (I-\widetilde{\mathcal{T}} \,,B_{R_0}\,, 0)$ is well defined and is equal to one.
Indeed, set $H_t (u):=t \,\widetilde{\mathcal{T}} u$ for all $t \in [0,1]$ and notice that
$(I-H_t)\,u=0 \;\Leftrightarrow\; u=t\,  \widetilde{\mathcal{T}} u \in B_{R_0}$. Then $I-H_t$ does not vanish on $\partial B_{R_0}\,$ and
\begin{align*}
\mathrm{deg} (I-\widetilde{\mathcal{T}} \,,B_{R_0}\,, 0)=\mathrm{deg} (I-H_1\,,B_{R_0}\,, 0)=\mathrm{deg} (I-H_0\,,B_{R_0}\,, 0)=\mathrm{deg} (I,B_{R_0}\,, 0)=1\, .
\end{align*}

Therefore, $\widetilde{\mathcal{T}}$ has a fixed point $u\in E$, which is a solution of $(\widetilde{\mathrm{P}})$. By claim \ref{claim1}, the first existence statement in theorem \ref{th2.1} is proved in the case $\kappa=\iota=1$.

If $\alpha$ and $\beta$ are in the general case as the maximum and minimum of strong sub and supersolutions, respectively, we define $\widetilde{f}$ as
\begin{align*}
\widetilde{f}(x,r,\vec{p}):=
\begin{cases}
\max_{1\leq i\leq \kappa}\overline{f}(x,\alpha_i (x),D\alpha_i (x))\, ,&\mathrm{if}\;\; r<\alpha (x) \\
\;\;\overline{f}(x,r,\vec{p}) \, ,\, &\mathrm{if}\;\;\; \alpha\leq r \leq \beta (x) \\
\min_{1\leq j\leq \iota}\overline{f}(x,\beta_j (x),D\beta_j (x))\, , &\mathrm{if}\;\; r>\beta (x)
\end{cases}
\end{align*}
and consider $R>\max \{ \|\alpha_i\|_{E},\|\beta_j\|_{E}\,; \; 1\leq i\leq \kappa, \, 1\leq j\leq \iota \}$. In claim \ref{claim1}, choose $i\in \{1,\ldots,\kappa\}$ such that $\min_{\overline{\Omega}}\, (u-\alpha)=(u-\alpha_i)(x_0)$, define $v=u-\alpha_i$ and replace $\alpha$ by $\alpha_i\in W^{2,p}(\Omega)$ until the end of the proof, observing that $x\in \widetilde{\Omega}=\{v<0\}$ implies $u(x)<\alpha (x)$. The rest of the proof is exactly the same.

\vspace{0.5cm}

\textit{Part 2. Degree computation in $\mathcal{S}$ under strictness of  $\alpha , \beta$ -- proof of Theorem \ref{th2.1} $(i)$.}

Suppose $\alpha , \beta$ are strict, and consider the set $\mathcal{S}$ as in the statement of theorem \ref{th2.1}.
Since there exists a solution $u \in C_0^1 (\overline{\Omega})$ of (P) with $\alpha\leq u \leq\beta$ in $\Omega$, by definition \ref{def2.1} we have $\alpha\ll u \ll \beta$ in $\Omega$ and so $\mathcal{S}$ is a nonempty open set in $C_0^1 (\overline{\Omega})$. Further, from part 1, we see that all fixed points of $\widetilde{\mathcal{T}}$ are in $\mathcal{S}\subset B_R \subset B_{R_0}\,$, the degree over $B_{R_0}$ is equal to $1$ and solutions of (P) and ($\widetilde{\mathrm{P}}$) in $\mathcal{S}$ coincide, leading to
\begin{align*}
\mathrm{deg} (I-\mathcal{T},\mathcal{S}, 0)=\mathrm{deg} (I-\widetilde{\mathcal{T}},\,\mathcal{S}, 0)=\mathrm{deg} (I-\widetilde{\mathcal{T}},\,B_{R_0}\,, 0)=1.
\end{align*}

\vspace{0.1 cm}
\textit{Part 3. Existence of extremal solutions under \eqref{Hstrong} -- proof of $(ii)$. }

Denote by $\mathcal{H}$ the set of fixed points of $\mathcal{T}$ belonging to the order interval $[\alpha,\beta]$. In claim \ref{claim1} we saw that this set contains the set of fixed points of $\widetilde{\mathcal{T}}$. The converse is also true since any solution of $(\mathrm{P})$ in the order interval $[\alpha,\beta]$ satisfies $\overline{f}(x,u,D u)=f(x,u,D u)$ and \eqref{u<R}, hence is a solution of $(\widetilde{\mathrm{P}})$, i.e. we have
$$
\mathcal{H}=\{ u \in E\, ; \; u=\mathcal{T}u\, , \; \alpha\leq u \leq\beta\;\; \mathrm{in}\;\Omega\}=\{u \in E\, ; \; u=\widetilde{\mathcal{T}} u \}\, .
$$
This set is nonempty by part 1. Also, from that proof, we know that for any  sequence $(u_k)_{k\in \n}\subset E$ with $u_k\rightarrow u$ in $E$ there exists $U\in E$ such that $\widetilde{\mathcal{T}}u_k \rightarrow U=\widetilde{\mathcal{T}} u$ in $E$. Then, if additionally $u_k=\widetilde{\mathcal{T}} u_k\,$, we obtain that $u=\widetilde{\mathcal{T}}u$ thus $\mathcal{H}$ is a compact set in $E$.

Now consider, for each ${u\in\mathcal{H}}\,$, the set $C_u:=\{ z\in\mathcal{H}\, ;\;\; z \leq u \;\;\mathrm{in}\;\Omega \}$.

\begin{claim} \label{claim2}
The family $\{ C_u \}_{u\in\mathcal{H}}\,$ has the finite intersection property, i.e the intersection of any finite number of sets  $\mathcal{C}_{u}$ is not empty.
\end{claim}

\begin{proof}
Let $u_1 , ... ,u_\kappa \in\mathcal{H}$. Observe that $\widetilde{\beta}:=\min_{1\leq i \leq \kappa} \, u_i$ is an $L^p$-viscosity supersolution of (P), with $\alpha\leq\widetilde{\beta} \leq \beta$ in $\Omega\,$. Furthermore, under hypothesis \eqref{Hstrong}, such $\widetilde{\beta}$ is a minimum of strong solutions of (P), which is exactly what we need to use part 1 of the above proof to obtain the existence of a solution $v$ of (P) with $\alpha\leq v \leq\widetilde{\beta} \leq\beta$ in $\Omega$, i.e. $v\in\mathcal{H}$ and $v\leq u_i\,$, for every $\, i \in \{1,...,\kappa\}$, so $v \in \cap_{1\leq i \leq \kappa} \,\mathcal{C}_{u_i}\neq \emptyset$.
\qedhere{\textit{Claim \ref{claim2}. }}
\end{proof}

By the definition of compacity of $\mathcal{H}$ by open covers and claim \ref{claim2}, there exists $\underline{u} \in \cap_{u\in\mathcal{H}}\,\mathcal{C}_u\,$ (see, for example, theorem 26.9 in \cite{Munkres}). But then there exists a solution $\underline{u}$ of (P) with $\alpha \leq\underline{u} \leq u$ in $\Omega$, for all $u\in\mathcal{H}$. Analogously we prove the existence of $\overline{u}$, with $u\leq \underline{u}\leq\beta$ in $\Omega$, for every $ u\in\mathcal{H}$.
 \qedhere{\,\textit{Theorem \ref{th2.1}.}}
\end{proof}

\begin{rmk}\label{remark th 2.1 com c(x,u)}
The conclusion of the theorem \ref{th2.1} is still true if, instead of $c(x)$, we have some $c(x,u)$ such that $c(x,u)u=c(x)T_a u$, where $T_a$ is a truncation of $u$, i.e. $T_a(u)=u$ for $u\geq a$, $T_a (u)=a$ for $u<a$. In this case, $|T_a u_k-T_a u|\leq |u_k -u|$ and the rest of the proof carries on in the same way.
\end{rmk}


\section{A priori bounds}\label{secao a priori bounds}

In this section we prove Theorem \ref{apriori}.
Now we look at our family of problems \ref{Plambda}
for $\lambda >0$, assuming $c\gneqq 0$ and that the matrix $M$ satisfies \ref{M}. Here we consider \ref{SC0} i.e we suppose that all coefficients of the problem \ref{Plambda} are bounded and $d\equiv 0$. With the latter the zero order term in \ref{Plambda} is explicit, so we can obtain a clear behavior of the solutions with respect to $\lambda $.

\begin{teo} \label{aprioriLinfty}
Let $\Omega\in C^{1,1}$ be a bounded domain. Suppose \ref{SC0}, \ref{H0} hold and let $\Lambda_1, \, \Lambda_2$ with $0<\Lambda_1<\Lambda_2$. Then every $L^p$-viscosity solution  $u$ of \ref{Plambda} satisfies
$$
\|u\|_\infty \leq C\, , \;\textrm{ for all } \lambda\in [\Lambda_1,\Lambda_2],
$$
where $C$ depends on $n,p,\mu_1, \Omega, \Lambda_1,\Lambda_2,  \|b\|_{L^\infty(\Omega)},\|c\|_{L^\infty(\Omega)}, \|h\|_{L^\infty(\Omega)},\|u_0\|_{L^\infty(\Omega)},\lambda_P,\Lambda_P$, the $C^{1,1}$ diffeomorphism that describes the boundary, and  the set where $c>0$.
\end{teo}


The proof of theorem \ref{aprioriLinfty} uses and develops  the  ideas sketched in \cite{note}, adding some improvements in order to remove restrictions on the size of $c$.
We start by proving that all supersolutions stay uniformly bounded from below, even when $\lambda$ is close to zero.

\begin{prop} \label{aprioribelowLinfty}
Let $\Omega\in C^{1,1}$ be a bounded domain. Suppose \ref{SC0} and let $\Lambda_2>0$. Then every $L^p$-viscosity supersolution $u$ of \ref{Plambda} satisfies
$$
\|u^-\|_\infty \leq C\, , \; \textrm{ for all } \lambda\in [0,\Lambda_2]
$$
where $C$ depends only on $n, p,\mu_1, \Omega,\Lambda_2,  \|b\|_{L^\infty(\Omega)},\|c\|_{L^\infty(\Omega)},  \|h^-\|_{L^\infty(\Omega)}$, $\lambda_P,\Lambda_P$.
\end{prop}

\begin{proof}
First observe that both $-u$ and $0$ are $L^p$-viscosity subsolutions of
\begin{align*}
-F(x,U,DU,D^2U)\leq \lambda c(x)U-\mu_1 |D U|^2 +h^-(x) \;\; \; \textrm{in} \;\;\; \Omega.
\end{align*}
Then, using \ref{SC0}, both are also $L^p$-viscosity subsolutions of
\begin{align*}
\left\{
\begin{array}{rclcc}
\mathcal{M}^+ (D^2 U)+b(x)|DU|-\mu_1 |D U|^2 &\geq & -\lambda c(x)U-\-h^-(x) &\mbox{in} & \Omega\\
U &\leq & 0 &\mbox{on} &\partial\Omega
\end{array}
\right.
\end{align*}
and so is $U:=u^-=\max \{-u,0\}$, as the maximum of subsolutions. Moreover, $U\geq 0$ in $\Omega$ and $U=0$ on $\partial\Omega$. We make the following exponential change
$$
w:=\frac{1-e^{-mU}}{m}, \;\;\; \mathrm{with}\;\; m=\frac{\mu_1}{\Lambda_P}
$$
where $\Lambda_P$ is the constant from the definition of Pucci's operators. From lemma \ref{lemma2.3arma}, we know that $w$ is an $L^p$-viscosity solution of
\begin{align}\label{Qlambda}
\tag*{($Q_\lambda$)}
\left\{
\begin{array}{rclcc}
-\mathcal{L}_1^+[w] &\leq & h^-(x)+\frac{\lambda}{m} c(x)\,|\mathrm{ln}(1-mw)|(1-mw)  &\mbox{in} & \Omega \\
w &=& 0 &\mbox{on} & \partial\Omega
\end{array}
\right.
\end{align}
where $\mathcal{L}_1^+[w]:=\mathcal{L}^+[w]-mh^-(x)w$. Notice that the logarithm above is well defined, since
$$
0\leq \,w=\frac{1}{m}\{1-e^{-mU}\}\,< \frac{1}{m} \;\;\;\; \mathrm{in}\;\; \Omega.
$$

Now set $w_1:=\frac{1}{m}(1-e^{-mu_1^-})$, where $u_1$ is some fixed supersolution of \ref{Plambda}, $\lambda\geq 0$ (if there was not such supersolution, we have nothing to prove). Then, by the above, $w_1\in [0,1/m)$ is an $L^p$-viscosity solution of ($Q_\lambda$). Define
$$
\overline{w}:=\sup \mathcal{A},\;\;\mbox{where}\;\; \mathcal{A}:=\{\,w: \;w\textrm{ is an } L^p \textrm{-visc.\ solution of (}Q_\lambda\mathrm{)};\; 0\leq w < 1/m \; \textrm{ in } \Omega  \,\}.
$$
Then $\mathcal{A}\neq \emptyset$ since $w_1 \in \mathcal{A}$, and $w_1\leq \overline{w}\leq 1/m$ in $\Omega$. Also, as a supremum of subsolutions (locally bounded, since it belongs to the interval $[0,1/m]$), $\overline{w}$ is an $L^p$-viscosity solution of the first inequality in ($Q_\lambda$). Clearly, $\overline{w}=0$ on $\partial\Omega$.

Observe that $\mathcal{L}_1^+$ is a coercive operator and the function
$$f(x):=f_\lambda(x,\overline{w}(x))=h^-(x)+\frac{\lambda}{m} c(x)\,|\mathrm{ln}(1-m\overline{w})|(1-m\overline{w}) \,\in L^p_+(\Omega)$$
with
$$
\|f^+\|_{L^p(\Omega)}\leq \|h^-\|_{L^p(\Omega)} + \frac{\Lambda_2}{m} \|c\|_{L^p(\Omega)} \,C_0
$$
since $A(\overline{w}):=|\mathrm{ln}(1-m\overline{w})|\,(1-m\overline{w})\leq C_0$. Indeed, from $\lim_{t\rightarrow 0^+}{t\,\mathrm{ln}t}=0$ there exists a $\delta \in (0,1)$ such that $t|\mathrm{ln}t|\leq 1, $ for all $0<t<\delta$ i.e. $A(\overline{w})\leq 1$ when $x\in \left\{ \overline{w} > \frac{1-\delta}{m} \right\}$. If $x\in \left\{ \overline{w} \leq \frac{1-\delta}{m} \right\}$ then $A(\overline{w})\leq |\mathrm{ln} \delta |$ (notice that $1-m\overline{w} \leq 1$) so take $C_0= \max \{1, |\mathrm{ln} \delta | \}$.\medskip

Therefore, by the proof of the boundary Lipschitz bound (see theorem 2.3 in \cite{B2016}),
\begin{align*}
\overline{w} (x)\leq C \|f^+\|_{L^p(\Omega)} \,\mathrm{dist}(x,\partial\Omega) \rightarrow 0 \quad \mathrm{as} \;\;\;x\rightarrow\partial\Omega
\end{align*}
and so $\overline{w}\not\equiv \frac{1}{m}$. Observe that the function $\overline{w}$ can be equal to $1/m$ at some interior points.

If there were a sequence of supersolutions $u_k$ of ($P_\lambda$) in $\Omega$ with unbounded negative parts, then there would exist a subsequence such that
$$
 u_k^-(x_k)=\|u_k^-\|_\infty \xrightarrow[{k\rightarrow \infty}]{} +\infty, \;\; x_k\in \overline{\Omega}, \;\; x_k\xrightarrow[{k}]{} x_0 \in \overline{\Omega}
$$
with $x_k\in\Omega$ for large $k$, since $u_k\geq 0$ on $\partial\Omega$. Then the respective sequence
$$w_k(x_k)=\frac{1}{m} \{ 1-e^{-mu_k^-(x_k)} \} \xrightarrow[{k\rightarrow \infty}]{} \frac{1}{m}\,, \;\; w_k \in \mathcal{A} $$
i.e. for every $\varepsilon >0$, there exists some $k_0\in \n$ such that
$$
\frac{1}{m} \geq \overline{w}(x_k)\geq w_k(x_k)\geq \frac{1}{m}-\varepsilon , \;\; \textrm{ for all }\, k\geq k_0
$$
thus there exists $\lim_k \overline{w}_(x_k)=\frac{1}{m}$ and also
$$
\overline{w}(x_0)\geq \displaystyle\varliminf_{x_k\rightarrow x_0} \overline{w} (x_k)=\lim_{k\rightarrow\infty} \overline{w} (x_k)=\frac{1}{m}.
$$
Hence $x_0 \in\Omega$, since $\overline{w}=0$ on $\partial\Omega$, and $\overline{w}(x_0)=\frac{1}{m}$.

Finally, define $z:=1-m\overline{w}$. Then $z$ is an $L^p$-viscosity supersolution of
$$
-\mathcal{L}_1^- z \geq -\lambda c(x)\, |\mathrm{ln} z| z \quad \mathrm{in}\;\; \Omega, \qquad z\gvertneqq 0\mbox{ in }\Omega, \quad z(x_0)=0,
$$
where $\mathcal{L}_1^-:=\mathcal{L}^- -mh^-$ is a coercive operator. But this contradicts the following nonlinear version of the SMP.
\end{proof}

\begin{lem} \label{vazquez particular} 
Set $\mathcal{L}_1^-[u]:=\mathcal{M}^-(D^2u)-\gamma\,|Du| -du$, for a constant $d\geq 0$.
Let $f\in C[0,+\infty)$ be defined by $f(s)=a\, s\,|\mathrm{ln}s|$ if $s>0$, $f(0)=0$, where $a\geq 0$. Then, the SMP holds for the operator $\mathcal{L}_1^-[\cdot] - f(\cdot)$, i.e. if $u$ is a $C$-viscosity solution of
\begin{align*}
\left\{
\begin{array}{rclcc}
\mathcal{L}_1^- [u] &\leq & f(u) &\;\;\mbox{in}&\Omega \\
u &\geq & 0  &\;\;\mbox{in} &\Omega
\end{array}
\right.
\end{align*}
then either $u>0$ in $\Omega$ or $u\equiv 0$ in $\Omega$.
\end{lem}

This lemma can be seen as a form of the  Vazquez's strong maximum principle for our operators, since one over the square root of the primitive of $|\mathrm{ln}z|z$ is not integrable at $0$. Actually, it is not difficult to check that  we can also have a term $|\mathrm{ln}s|^b$ in  lemma \ref{vazquez particular},  for $b<2$. The proof of lemma \ref{vazquez particular} is given in the appendix.

Note that we  apply lemma  \ref{vazquez particular} with $d=m\,\|h^-(x)\|_{\infty}$ and $a=\Lambda_2\, \|c\|_\infty$.
\medskip

Before giving the proof of Theorem \ref{aprioriLinfty} we recall that the class of equations we study is invariant with respect to diffeomorphic changes of the spatial variable. In particular we can assume that the boundary of $\Omega$ is a hyperplane in a neighborhood of any given point of $\partial\Omega$. Indeed, straightening of the boundary leads to an equation of the same type, with bounds on the norms of the coefficients depending also on the $C^{1,1}$-norm of $\partial \Omega$.

\begin{proof}[\textit{Proof of Theorem \ref{aprioriLinfty}.}]
Fix $\Lambda_1, \Lambda_2$ with $0<\Lambda_1<\Lambda_2$. From proposition \ref{aprioribelowLinfty}, there exists a constant $C_1 >0$ such that
\begin{align} \label{u-leqC_1}
u^-\leq C_1 \, , \;\;\mathrm{for \;every\; supersolution\;} u\; \mathrm{of\; (}P_\lambda\mathrm{),}\; \textrm{ for all } \lambda\in [0,\Lambda_2].
\end{align}
Suppose then, in order to obtain a contradiction, that solutions are not bounded from above in $[\Lambda_1,\Lambda_2]$, by picking out a sequence $u_k$ of $L^p$-viscosity solutions of \ref{Plambda} such that
$$
 u_k^+(x_k)\xrightarrow[{k\rightarrow \infty}]{} +\infty, \;\; x_k\in \overline{\Omega}, \;\; x_k\xrightarrow[{k}]{} x_0 \in \overline{\Omega}.
$$
where $x_k$ is the point of maximum of $|u_k|$ in $\overline{\Omega}$, i.e. $\|u_k\|_\infty =|u(x_k)|$, $u_k^-(x_k)\in [0,C_1]$.

We claim that, up to changing the blow-up limit point $x_0$, we can suppose that there is a ball around $x_0$ in which $c$ is not identically zero.

To prove the claim, consider $G$, a maximal domain such that $c\equiv 0$ in $G$.
Obviously there is no need of such argument if $|\{c=0\}|=0$ or even if  $c\gneqq 0$ in a neighborhood of $x_0$. Suppose, hence, that $x_0$ is an interior point of $G$, and so $x_k\in G$ for large $k$ (considering a half ball in $G$ if $x_0\in \partial\Omega$, after a diffeomorphic change of independent variable which straightens the boundary).
Notice that both $u_k$ and $u_0$ satisfy the same equation in $G$, in the $L^ p$-viscosity sense, for each $k\in\n$.
Now, since there is no zero order terms of this equation in $G$, both $v_k:=u_k-\inf_{\partial G}u_k$ and $v_0:=u_0-\sup_{\Omega} u_0$ also satisfy this same equation, with $v_k\leq 0\leq v_0$ on $\partial G$, being, thus, $L^p$-viscosity sub and supersolutions of $(P_0)$, respectively, with $v_0$ strong. We apply lemma \ref{lema subsol strong always less supersol} to obtain that $v_k\leq v_0$ in $G$ and, in particular, for large $k$,
$$
 u_k^+ \geq u_k^+ (x_k) -2\|u_0\|_{L^\infty (\Omega)}-C_1 \xrightarrow[k\rightarrow\infty]{} +\infty \;\;\textrm{ on }\; {\partial G}.
$$
This means that there is blow-up also at the boundary of $G$, in the sense that there exists a sequence $y_k\in \partial G$ with $u_k^+(y_k)\rightarrow +\infty$ and $y_k\to y_0\in\partial G$, as $k\rightarrow \infty$.
Next, since $G$ is maximal, so $\partial G\subset \partial\Omega \cup \partial (\{c=0\})$, and using $u_k=0$ on $\partial\Omega$, we have $y_k\in \partial (\{c=0\})$. Therefore, we can take a ball $B_r(y_0)$ centered at $y_0$ (or a half ball if $y_0\in \partial\Omega$)  which, by enlarging $r$ if necessary,  becomes a neighborhood which meets the set $\{c>0\}$; in another words, such that $c\gneqq 0$ in $B_r(y_0)$.
Hence, up to changing $x_k$ and $x_0$ by $y_k$ and $y_0$, we can suppose that $c\gneqq 0$ in $B_r(x_0)$, or in a half ball if $x_0\in \partial\Omega$, after straightening the boundary around $x_0$.

Suppose we are in the more difficult case of a half ball. For simplicity, and up to rescaling, say $c\gneqq 0$ in $B_1^+=B_1^+(x_0)$, with our equation being defined in $B_2^+(x_0)\subset \Omega$.

We take the convention of assuming that the constant $C$ may change from line to line and depends on $n,\,p,\,\lambda_P,\,\Lambda_P,\,\Lambda_1, \,\Lambda_2$, $\mu_1$, $\|b\|_{L^p(\Omega)}$, $\|h\|_{L^p(\Omega)}$, $\|c\|_{L^p(\Omega)}$ and $C_1$. The constant $C_1$ is fixed in \eqref{u-leqC_1} with its dependence described in the statement of proposition~\ref{aprioribelowLinfty}.

Notice that, from \eqref{u-leqC_1}, for every $L^p$-viscosity solution $u$ of \ref{Plambda}, the function $v:=u+C_1$ is a nonnegative $L^p$-viscosity solution of
\begin{align*}
\mathcal{M}^-(D^2 v)-b(x)|Dv| &\leq F(x,v-C_1,Dv,D^2v)\\
&=-\lambda c(x)v-\langle M(x)D v, D v \rangle -h(x)+\lambda c(x)C_1 \\
& \leq -\lambda c(x)v - \mu_1 |D v|^2 +\widetilde{h}(x)
\end{align*}
where $\widetilde{h}(x):=h^-(x)+\Lambda_2 \,c(x) C_1 \geq 0\,$, by \ref{SC0}. Thus, by lemma \ref{lemma2.3arma}, the function
\begin{align} \label{v1}
v_1:=\frac{1}{m_1} \,\{ e^{m_1v}-1\}\, ,\;\,\mathrm{where} \;\; m_1=\frac{\mu_1}{\Lambda_P}
\end{align}
is a nonnegative $L^p$-viscosity supersolution of
\begin{align}\label{eq v1}
\mathcal{L}^-_1 [v_1]\,\leq\, f_1 (x) \quad \mathrm{in}\;\;B_2^+
\end{align}
where $\,\mathcal{L}^-_1 [v_1]:=\mathcal{M}^-(D^2 v_1)-b(x)|Dv_1|-m_1\,\widetilde{h} (x) \,v_1\;$ and $\; f_1(x):=-\frac{\lambda}{m_1} c(x)(1+m_1v_1)\,\\ \mathrm{ln} (1+m_1v_1)+\,\widetilde{h} (x)\in L^p (\Omega)$ since $v_1 \in L^\infty (\Omega)$.

Notice that, in the set $B_2^+\cap \{f_1\geq 0\}$, we have
$0\leq \frac{\lambda}{m_1}c(x)(1+m_1v_1)\,\mathrm{ln} (1+m_1v_1)\leq\widetilde{h}$
 and $f_1^+=|f_1|\leq \frac{\lambda}{m_1}c(x)(1+m_1v_1)\,\mathrm{ln} (1+m_1v_1) + \widetilde{h}\leq 2\,\widetilde{h}$, so
\begin{align}\label{eq f1+}
 \|f_1^+ \|_{L^p(B^+_2)}\leq 2 \|\widetilde{h} \|_{L^p(B^+_2)}\leq 2 \|h^- \|_{L^p(\Omega)}+2\Lambda_2\|c \|_{L^p(\Omega)} C_1\leq C.
\end{align}

Then, using proposition \ref{BQSMP} applied to \eqref{eq v1}, we obtain positive constants $ c_0, \, C_0$ and $\varepsilon\leq 1$, depending on $n,\,\lambda_,\,\Lambda_P, \,p$ and $\|b\|_{L^p(\Omega)}$, such that
\begin{align*}
I:\,=&\,\inf_{B_1^+} \frac{v_1}{x_n} \geq c_0 \left( \int_{B_{3/2}^+} (f_1^-)^\varepsilon \right)^{1/\varepsilon} -C_0 \,\|f_1^+ \|_{L^p(B^+_2)} \\
=& \; c_0 \left( \int_{B_{3/2}^+} \left\{ \left( \frac{\lambda}{m_1} c(x)(1+m_1v_1)\,\mathrm{ln} (1+m_1v_1) - \widetilde{h} (x) \right)^+\right\}^\varepsilon \,\right)^{1/\varepsilon} - C \\
\geq & \, c_0 \,\inf_{B_1^+} \frac{v_1}{x_n} \left( \int_{B^+_1} \left(\left( \lambda c(x) \frac{1+m_1v_1}{m_1 v_1} \,x_n \,\mathrm{ln} (1+m_1v_1) -  \widetilde{h} (x)\frac{1+m_1v_1}{ v_1} \,x_n \right)^+\right)^\varepsilon \right)^{\frac{1}{\varepsilon}} - C \\
\geq & \; c_0 \; I \left( \int_{B^+_1} \left\{ \left( \lambda c(x) \,x_n \,\mathrm{ln} (1+I\,m_1 x_n) -  m_1 \,\widetilde{h} (x) \,x_n \right)^+ \right\}^\varepsilon \right)^{1/\varepsilon} - C
\end{align*}
using \eqref{eq f1+} and that $v_1(x)\geq x_n I$ for all $x\in B^+_1\,$. Thus,
\begin{align} \label{I.leqC2}
I \,\left\{ c_0 \left( \int_{B_1^+} x_n^\varepsilon \,\left( \left( \lambda c(x) \,\mathrm{ln} (1+I\,m_1x_n) -  m_1 \,\widetilde{h} (x) \, \right)^+ \right)^\varepsilon \right)^{1/\varepsilon} -1 \right\} \leq \bar C.
\end{align}
We claim that this is only possible if $I \leq C$, with a constant that does not depend on $v_1$ (and consequently on $u$), neither $\lambda\in [\Lambda_1,\Lambda_2]$.
Indeed, if this was not the case, we would obtain a sequence of supersolutions $v_1^k$ of $\mathcal{L}^-_1 [v_1^k]\,\leq\, f_1^k (x)$ in $B_2^+$ such that $I_k:=\inf_{B_1^+} \frac{v_1^k}{x_n}\rightarrow +\infty$ when $k \rightarrow +\infty$ and \eqref{I.leqC2} holding with $I$ replaced by $I_k\,$. So, up to a subsequence and renumbering, we can assume that $I_k\geq k^2$ and $\frac{\bar C}{k^2}\leq 1$ for all $k\geq k_0$, from where we obtain that
\begin{align*}
\int_{B_1^+} x_n^\varepsilon \,\left( \left( \lambda c(x) \,\mathrm{ln} (1+I_k\,m_1 x_n)  -  m_1 \,\widetilde{h} (x) \, \right)^+ \right)^\varepsilon\,\leq \,c_0^{-\varepsilon }\left( 1+\dfrac{\bar{C}}{k^2}\right)^{\varepsilon}\leq C
\end{align*}
and finally, using $\lambda\geq\Lambda_1$,
\begin{align} \label{Ik.leqC2}
\int_{B_1^+\cap \{x_n\geq 1/k\}} x_n^\varepsilon \,\left( \left( \Lambda_1 c(x)  -  m_1 \,\frac{\widetilde{h} (x)}{\mathrm{ln}(1+m_1k)} \, \right)^+ \right)^\varepsilon\,\leq \,\dfrac{C}{\mathrm{ln}(1+m_1k)}.
\end{align}
Taking the limit when $k \rightarrow +\infty$ we have $\int_{B_1^+} (x_n\,c(x))^\varepsilon \d x =0$, since $\Lambda_1>0$, which contradicts $c(x)\gneqq 0$ in $B_1^+$. More precisely, for the limit in \eqref{Ik.leqC2} we can use, for example, the  dominated convergence theorem: for $\varepsilon=1$ this is obvious; for $0<\varepsilon<1$ we use Young's inequality to estimate
$$
\left( \left( \lambda c(x)  -  m_1 \,{\widetilde{h} (x)}/{\mathrm{ln}(1+m_1k)} \, \right)^+ \right)^\varepsilon \leq \left( \lambda c(x)  -  m_1 \,{\widetilde{h} (x)}/{\mathrm{ln}(1+m_1k)} \, \right)^+ + \,1
$$
ensuring the desired convergence.
In this way we have gotten the claim $\inf_{B_1^+} \frac{v_1}{x_n} \leq C.$ 

Thus, by theorem \ref{BWHI} applied to \eqref{eq v1}, we have that there exists other positive constants $\varepsilon, \,c_0, \, C_0$, depending on $n,\,\lambda_P,\,\Lambda_P, \,p$ and $\|b\|_{L^p(\Omega)}$, such that
\begin{align} \label{estim.v_1.epsilon}
 \left( \int_{B^+_{3/2}} {v_1}^\varepsilon \right)^{1/\varepsilon} \leq c_0 \left( \int_{B^+_{3/2}} \left(\frac{v_1}{x_n}\right) ^\varepsilon \right)^{1/\varepsilon} \leq  \;\inf_{B_1^+} \frac{v_1}{x_n} +C_0 \,\|f_1^+\|_{L^p(B^+_2)} \leq C .
\end{align}

Now we go back to $u$ and define
\begin{align} \label{v2}
v_2:=\frac{1}{m_2} \,\{ e^{m_2 u}-1\}\, ,\;\,\mathrm{with} \;\; m_2=\frac{\mu_2}{\lambda_P}
\end{align}
which by lemma \ref{lemma2.3arma} and \ref{SC0} is an $L^p$-viscosity solution of
\begin{align}\label{eq v2}
\left\{
\begin{array}{rclcc}
\mathcal{M}^+(D^2 v_2)+b(x)|Dv_2|+\nu (x)v_2 &\geq & -h^+(x) & \mbox{in} & B_2^+ \\
v_2 &=& 0  & \mbox{on} & B_2^0
\end{array}
\right.
\end{align}
where $\nu (x):=\frac{\lambda}{m_2v_2}c(x)(1+m_2v_2)\, \mathrm{ln}(1+m_2v_2)\in L^p (B_2^+)$. \medskip

Notice that by the definitions \eqref{v1} and \eqref{v2}
\begin{align}\label{relacao.v1ev2}
v_2=\frac{1}{m_2} \left\{\,(1+m_1v_1)^{\frac{m_2}{m_1}} \,e^{-m_2C_1} -1 \,\right\}.
\end{align}

As in \cite{CJ}, observe that
\begin{align}\label{estim.ln}
\lambda c(x) \left| \frac{(1+m_2v_2)}{m_2v_2}\, \mathrm{ln} (1+m_2v_2) \right| \leq C_s \, c(x) \left( 1+ |v_2|^s \right)
\end{align}
for any $s>0$. Now, if we take $s=\varepsilon\, \frac{m_1}{m_2}\frac{p-n}{p(p+n)}$ and $p_1= \frac{p+n}{2}\in (n,p)$ then, by Holder's inequality, the right hand side in \eqref{estim.ln} belongs to $L^{p_1} (B_2^+)$ and

\begin{align} \label{escolha de r}
\|\,c\, |v_2|^s \|_{L^{p_1}(B_{2}^+)} &\leq \|c \|_{L^{p}(B_{2}^+)} \|\, |v_2|^r \|_{L^{p_2}(B_{2}^+)} \, , \quad \mathrm{with}\;\;\; \frac{1}{p_1}=\frac{1}{p}+\frac{1}{p_2} \nonumber \\
&\leq \|c \|_{L^{p}(B_{2}^+)} \left( \int_{B_{2}^+} |v_2|^{\varepsilon \frac{m_1}{m_2}}\right)^{\frac{p-n}{p(p+n)}} \leq \, C
\end{align}
and then
\begin{align}\label{cota nu}
\|\nu \|_{L^{p_1}(B_2^+)}\leq C_{n,p}\,\|c\|_{L^p(\Omega)}+\|\,c\, |v_2|^s \|_{L^{p_1}(B_{2}^+)} \leq C.
\end{align}

Now, the uniform bound on \eqref{cota nu} allows us to use proposition \ref{BLMP} applied to \eqref{eq v2}, for $v_2$ and $r=\varepsilon \frac{m_1}{m_2}$ as in  \eqref{escolha de r}, to obtain
\begin{align*}
{v_2}^+ \leq C \, \left\{  \left( \int_{B_{3/2}^+} |v_2|^{r}\right)^{{1}/{r}} +\|h^+ \|_{L^p (\Omega)}  \right\} \leq\, C\;\;\textrm{in } B_1^+.
\end{align*}
Hence, $u^+=\frac{1}{m_2} \,\mathrm{ln} (1+m_2{v_2}^+)$ is uniformly bounded in $B_1^+$, for every $L^p$-viscosity solution of \ref{Plambda}, for all $\lambda\in [\Lambda_1,\Lambda_2]$.
\end{proof}

\section{Proofs of the main theorems}\label{secao proofs}

In  this section we  assume that  $c\gneqq 0$ and \ref{M} holds.

\subsection{Some auxiliary results}\label{subsection auxiliary results}

We start by constructing an auxiliary problem \eqref{Plambda,k}, for which we can assure that there are no solutions for large $k$, and such that $(P_{\lambda,0})$ reduces to the problem \ref{Plambda}. This is a typical  but essential (see \cite{BR}, \cite{CJ}) argument that allows us to find a second solution via degree theory, by homotopy invariance in $k$.

Fix $\Lambda_2>0$. Recall that proposition \ref{aprioribelowLinfty} gives us an a priori lower uniform bound $C_0$, depending only on $n,p,\lambda_P,\Lambda_P, \,\mu_1, \,\Omega,\,\Lambda_2, \, \|b\|_{L^\infty(\Omega)},\,\|c\|_{L^\infty(\Omega)}$ and $\|h^-\|_{L^\infty(\Omega)}$, such that
\begin{align} \label{u geq -C_0 for Plambda,k}
u\geq -C_0\,,\;\;\mathrm{for\; every\;} L^p \textrm{-viscosity  supersolution } u \; \mathrm{of} \; (P_\lambda),\; \textrm{ for all }\,\lambda\in [0,\Lambda_2].
\end{align}

Consider, thus, the problem
\begin{align} \tag{$P_{\lambda ,k}$} \label{Plambda,k}
\left\{
\begin{array}{rclcc}
-F[u]&=&\lambda c(x)u+h(x)+\langle M(x)D u,D u\rangle +k\,\widetilde{c}(x) &\mbox{in} & \Omega \\
u &=& 0 &\mbox{on} & \partial\Omega
\end{array}
\right.
\end{align}
for $k\geq 0$, $\lambda\in [0,\Lambda_2]$, $F$ satisfying \ref{SC0}, \eqref{Hbeta} and \eqref{ExistUnic M bem definido}, $M$ satisfying \ref{M}, $c\gvertneqq 0$, $c,h\in L^\infty(\Omega)$ and $\widetilde{c}\,$ being defined as
\begin{align} \label{ctilde e A for k=1}
\widetilde{c}(x)=\widetilde{c}_{\Lambda_2} (x):= A c(x)+h^-(x)+\Lambda_2\, C_0\, c(x) \in L^\infty_+ (\Omega),
\end{align}
with $A:={\lambda_1}/{m}\,$, $m={\mu_1}/{\Lambda_P}\,$, where $\lambda_1=\lambda_1^+\left( \mathcal{L}^-(c),\Omega \right)>0$ is the first eigenvalue with weight $c$ of the proper operator $\mathcal{L}^-$, associated to the positive eigenfunction $\varphi_1 \in W^{2,p}(\Omega)$, given by proposition  \ref{exist eig for L-c geq 0}.

Note that every $L^p$-viscosity solution of $(P_{\lambda,k})$ is also an $L^p$-viscosity supersolution of \ref{Plambda}, since $k\,\widetilde{c}\geq 0$, and so satisfies \eqref{u geq -C_0 for Plambda,k}.
From this and \eqref{ctilde e A for k=1} we have, for all $ k\geq 1$,
\begin{align} \label{ctilde e A for all k natural}
\lambda c(x)u+ h(x)+k\,\widetilde{c}(x)\geq -\Lambda_2C_0\,c(x)-h^-(x)+\widetilde{c}(x)= A c(x) \gneqq 0 \;\textrm{ a.e. in } \Omega .
\end{align}

\begin{lem} \label{lema Plambda,k has no solutions}
For each fixed $\Lambda_2>0$, $(P_{\lambda ,k})$ has no solutions for all $k\geq 1$ and $\lambda\in [0,\Lambda_2]$.
\end{lem}

\begin{proof}
First observe that every $L^p$-viscosity solution of $(P_{\lambda,k})$, for $\lambda\in [0,\Lambda_2]$, is positive in $\Omega$.
Indeed, from \eqref{ctilde e A for all k natural}, \ref{SC0} and $M\geq 0$, we have that $u$ is an $L^p$-viscosity solution of
\begin{align*}
\left\{
\begin{array}{rclcc}
\mathcal{L}^-[u] &\lneqq & 0 &\;\;\mbox{in}&\Omega \\
u &=& 0 &\;\;\mbox{on} & \partial\Omega,
\end{array}
\right.
\end{align*}
and this implies that $u\geq 0$ in $\Omega$ by ABP.
Then $u>0$ in $\Omega$ by SMP.

Assume, in order to obtain a contradiction, that $(P_{\lambda,k})$ has a solution $u$. Again by \ref{M}, \eqref{ctilde e A for all k natural} and \ref{SC0}, we see that $u$ is also an $L^p$-viscosity solution of
\begin{align*}
\left\{
\begin{array}{rclcc}
\mathcal{L}^-[u] &\leq & -\mu_1 |D u|^2 -Ac(x) &\mbox{in}& \Omega \\
u &>& 0 &\mbox{in} & \Omega,
\end{array}
\right.
\end{align*}
and from lemma \ref{lemma2.3arma},
\begin{align}\label{eq v L-(c) to absurd}
\left\{
\begin{array}{rclcc}
(\mathcal{L}^-+\lambda_1c) [v] &\leq & -Ac(x) &\mbox{in}& \Omega \\
v &>& 0  &\mbox{in} &\Omega,
\end{array}
\right.
\end{align}
using $mA=\lambda_1\,$, where $v=\frac{1}{m} \left\{ e^{mu} -1 \right\}$, for $m$ and $A$ from \eqref{ctilde e A for k=1}. Then \eqref{eq v L-(c) to absurd} and \eqref{eq exist eigen L-c}, together with proposition \ref{th4.1 QB}, yields $v=t\varphi_1$ for $t>0$. But this contradicts the first line in \eqref{eq v L-(c) to absurd}, since $Ac(x)\ngeq 0$ and $(\mathcal{L}^-+\lambda_1c) [t\varphi_1]=t(\mathcal{L}^-+\lambda_1c) [\varphi_1]=0$\, in $\Omega$.

\qedhere{\textit{\,Lemma \ref{lema Plambda,k has no solutions}.}}
\end{proof}

When we are assuming hypothesis \eqref{Hstrong} we just say solutions to mean strong solutions of $(\overline{P}_\lambda)$. However, sub and supersolutions of such equations, in general, are not assumed  strong (since we are considering the problem in the $L^p$-viscosity sense), unless specified. In order to avoid possible confusion, we always make explicit the notion of sub/supersolution we are referring to.

The next result is  important in degree arguments, bearing in mind the set $\mathcal{S}$ in theorem~\ref{th2.1}. It will play the role of the strong subsolution $\alpha$ in that theorem.

\begin{lem}\label{lemma 4.2}
Suppose \ref{SC0}, \eqref{Hbeta} and \eqref{Hstrong}. Then, for every $\lambda>0$, there exists a strong strict subsolution  $\alpha_\lambda$ of \ref{Plambda} which is strong minimal, in the sense that every strong supersolution $\beta$ of \ref{Plambda} satisfies $\alpha_\lambda\leq\beta$ in $\Omega$.
\end{lem}

When $u_0$ has a sign, we will see in the proofs of theorems \ref{th1.5} and \ref{th1.4} that $u_0$ can be taken as $\beta$ and $\alpha$, respectively, in theorem \ref{th2.1} for the problem \ref{Plambda}, for all $\lambda>0$.

\begin{proof}
Let $K$ be the positive constant from proposition \ref{aprioribelowLinfty} such that every $L^p$-viscosity supersolution $\beta$ of
\begin{align}\tag{$\widetilde{P}_\lambda$}
\left\{
\begin{array}{rclcc}
-F[\beta] &\geq & \lambda c(x)\beta+\langle M(x)D \beta,D \beta\rangle -h^-(x)-1 & \mbox{in}&\Omega \\
\beta &\geq & 0 &\mbox{on}& \partial\Omega
\end{array}
\right.
\end{align}
satisfies $\beta\geq - K$ in $\Omega$.
Let $\alpha_0$ be the strong solution of the problem
\begin{align}\label{eq alpha0}
\left\{
\begin{array}{rclcc}
\mathcal{L}^-[\alpha_0]&=& \lambda Kc(x)+h^-(x)+1 & \mbox{in}& \Omega \\
\alpha_0&=&0 &\mbox{on}& \partial\Omega,
\end{array}
\right.
\end{align}
given, for example, by proposition 2.4 in \cite{KSexist2009}.
Then, as the right hand side of \eqref{eq alpha0} is positive, by ABP, SMP and Hopf, we have $\alpha_0\ll 0$ in $\Omega$.

\begin{claim}\label{claim 1 lemma 4.2}
Every $L^p$-viscosity supersolution $\beta$ of \ref{Plambda} satisfies $\beta\geq \alpha_0$ in $\Omega$.
\end{claim}

\begin{proof}
First notice that $\beta$ is an $L^p$-viscosity supersolution of $(\widetilde{P}_\lambda)$ and so satisfies $\beta\geq -K$. Second, by \ref{SC0} and $M\geq 0$, $\beta$ is also an $L^p$-viscosity supersolution of
$$
-\mathcal{L}^-[\beta]\geq \lambda c(x)\beta+h(x)\geq -\lambda K c(x)-h^-(x) -1\quad\mathrm{in}\;\;\Omega
$$
and setting $v:=\beta-\alpha_0$ in $\Omega$, $v$ is an $L^p$-viscosity solution of
$$
\mathcal{L}^-[v]\leq \mathcal{M}^-(D^2\beta)+\mathcal{M}^+(-D^2\alpha_0)-b(x)|D\beta|+b(x)|D\alpha_0|=\mathcal{L}^-[\beta]-\mathcal{L}^-[\alpha_0]\leq 0
$$
since $\alpha_0$ is strong. Further, $v\geq 0$ on $\partial\Omega$ then, by ABP, $v\geq 0$ in $\Omega$.
\qedhere{\,\textit{Claim \ref{claim 1 lemma 4.2}.}}
\end{proof}

Set
\begin{align*}
\overline{c}\,(x,t)=
\begin{cases}
\;\;\;\;c(x) \quad &\mathrm{if}\; \;\;t\geq -K \\
-{K}\,c(x)/t \quad & \mathrm{if} \;\;\; t<-K.
\end{cases}
\end{align*}
Observe that
$0\leq \overline{c}\,(x,t)\leq c(x)\,$ a.e. in $\Omega$ and $\overline{c}(x,t)t\geq -Kc(x)$ for all $t\in\real$. Then,
$$
-F[\alpha_0]\leq -\mathcal{L}^-[\alpha_0]=-\lambda K c(x)-h^-(x)-1\leq \lambda \,\overline{c}\,(x,\alpha_0)\alpha_0+\langle M(x)D\alpha_0,D\alpha_0\rangle -h^-(x)-1
$$
since $M\geq 0$ and $\alpha_0$ is a strong subsolution of $(\overline{P}_\lambda)$.

 Consider the problem $(\overline{P}_\lambda)$, which we define as  the problem \ref{Plambda} with $c,h$ replaced by $\overline{c}=\overline{c}(x,u)$, $\overline{h}=-h^--1$.

Observe that we are in the situation of remark \ref{remark th 2.1 com c(x,u)}, since $\bar{c}(x,u)=c(x)T_{-K}u$ there, which allows us to use theorem \ref{th2.1} in order to obtain solutions of $(\overline{P}_\lambda)$.

Let $\beta_0$ be some fixed strong supersolution of \ref{Plambda} (if there were not strong supersolutions of \ref{Plambda}, the proof is finished). Then, by claim \ref{claim 1 lemma 4.2}, we have $\alpha_0\leq\beta_0$ in $\Omega$. Also, in that proof we observed that $\beta_0\geq - K$, then $\overline{c}\,(x,\beta_0)\equiv c(x)$ a.e. $x\in\Omega$, which means that
$$
-F[\beta_0]\geq \lambda c(x)\beta_0+\langle M(x)D\beta_0,D\beta_0\rangle +h(x)\geq \lambda\, \overline{c}\,(x)\beta_0+\langle M(x)D\beta_0,D\beta_0\rangle -h^- (x)-1.
$$
and so $\beta_0$ is a strong supersolution of $(\overline{P}_\lambda)$. By theorem \ref{th2.1} and remark \ref{remark th 2.1 com c(x,u)}, we obtain an $L^p$-viscosity solution $w$ of this problem, with $\alpha_0\leq w \leq \beta_0$ in $\Omega$, which is strong and can be chosen as the minimal solution in the order interval $[\alpha_0,\beta_0]$, by hypothesis \eqref{Hstrong}.

\begin{rmk}\label{obs w Plambda tilde}
Notice that $\overline{c}\,(x,t)t\geq c(x)t\,$ a.e. $x\in\Omega\,$ for all $t\in\real$, so
\begin{align}
-F[w]&= \lambda\, \overline{c} \,(x,w)w+\langle M(x)D w ,D w \rangle -h^-(x)-1\nonumber\\
&\geq \lambda \,c(x)w+\langle M(x)D w ,D w \rangle -h^-(x)-1,\label{loc1}
\end{align}
a.e. in $\Omega$, i.e. $w$ is also a strong supersolution of $(\widetilde{P}_\lambda)$.
\end{rmk}

\begin{claim} \label{claim 2 lemma 4.2}
For every $\beta$ strong supersolution of \ref{Plambda}, $\beta\geq w$ in $\Omega$.
\end{claim}

\begin{proof}
Let $\beta$ be any strong supersolution of \ref{Plambda}. As in the argument above for $\beta_0$, we have that $\beta$ is also a strong supersolution of $(\overline{P}_\lambda)$. Suppose that the conclusion is not verified, i.e. that there exists $x_0\in \Omega$ such that $\beta (x_0)<w (x_0)$ and define
$$
\beta_1:=\min_\Omega\, \{w\,,\beta\} \not\equiv w.
$$
Then $\beta_1$ is the minimum of strong supersolutions, hence itself is an $L^p$-viscosity supersolution of $(\overline{P}_\lambda)$ and of $(\widetilde{P}_\lambda)$, by remark \ref{obs w Plambda tilde}. Following the same lines as in claim \ref{claim 1 lemma 4.2}, $\beta_1\geq \alpha_0$ in $\Omega$. Thus, by theorem \ref{th2.1} and remark \ref{remark th 2.1 com c(x,u)}, there exists an $L^p$-viscosity solution $w_1$ of $(\overline{P}_\lambda)$, strong by \eqref{Hstrong}, such that
$\alpha_0\leq w_1\leq \beta_1 \lneqq w\leq \beta_0$ in $\Omega$,
which gives a contradiction with the minimality of $w$.
\qedhere{\,\textit{Claim \ref{claim 2 lemma 4.2}.}}
\end{proof}

\begin{claim} \label{claim 3,4 lemma 4.2}
$w$ is a strong strict subsolution of \ref{Plambda}.
\end{claim}

\begin{proof}
From remark \ref{obs w Plambda tilde}, $w\geq -K$ and $\overline{c}\,(x,w)\equiv c(x)$, which implies that $w$ actually satisfies \eqref{loc1} with  equality, from where
\begin{align}
-F[w]< \lambda c(x)w+\langle M(x)D w ,D w \rangle +h(x) \;\;\;\mathrm{a.e.\;in}\;\Omega,
\end{align}
with $w=0$ on $\partial\Omega$. Thus $w$ is a strong subsolution of \ref{Plambda}. What remains to be proved is that $w$ is strict, in order to choose $\alpha_\lambda$ as $w$.
Therefore, in the sense of definition \ref{def2.1}, let $u\in E$ be an $L^p$-viscosity supersolution of \ref{Plambda} with $u\geq w$ in $\Omega$. Then, since $w$ is strong, $U:=u-w\geq 0$ in $\Omega$ is an $L^p$-viscosity supersolution of
\begin{align*}
-\mathcal{L}^-[U]\geq \lambda c(x)U+\langle M(x)D U,D U\rangle +\langle M(x)D U,D w\rangle +\langle M(x)D w,D U\rangle
\end{align*}
in $\Omega$, using \ref{SC0} and $0\leq M(x)\leq \mu_2$. Hence, for $\widetilde{b}=b+2\mu_2 |Dw|\in L^p_+(\Omega)$, $U$ satisfies
\begin{align*}
\left\{
\begin{array}{rclcc}
\mathcal{M}^-(D^2U)-\widetilde{b}(x)|DU|&< & 0 & \mbox{in} & \Omega \\
U &\geq & 0 & \mbox{in}& \overline{\Omega}
\end{array}
\right.
\end{align*}
in the $L^p$-viscosity sense and by SMP, $U>0$ in $\Omega$. If there exists $x_0\in\partial\Omega$ with $U(x_0)=0$, Hopf lemma implies that $\partial_\nu U(x_0)>0$. Then, $U\gg 0$ in $\Omega$.
\qedhere{\,\textit{Claim \ref{claim 3,4 lemma 4.2}.}}
\end{proof}
\qedhere{\,\textit{Lemma \ref{lemma 4.2}.}}
\end{proof}

\subsection{Proof of Theorem \ref{th1.1,1.2,1.3}}

Suppose, at first, \eqref{SC}, \eqref{Hbeta} and \eqref{ExistUnic M bem definido}.

We start proving the first statement in theorem \ref{th1.1,1.2,1.3}, about existence of solutions for $\lambda<0$.
Set $\alpha:=u_0-\|u_0\|_{\infty}$ and $\beta:=u_0+\|u_0\|_{\infty}$. Thus, $\alpha,\beta$ is a pair of strong sub and supersolutions of \ref{Plambda}, for each $\lambda<0$, with $\alpha\leq\beta$ in $\Omega$. Indeed, using \eqref{SC}, $\alpha\leq u_0$ and $0\leq\lambda c(x)\alpha$, we have
\begin{align} \label{eq1 teo1}
\left\{
\begin{array}{rclcc}
-F[\alpha]&\leq & -F[u_0]\;\leq\; \lambda c(x)\alpha+h(x)+\langle M(x)D \alpha,D\alpha\rangle &\mbox{in}& \Omega\\
\alpha &\leq& u_0\;=\;0 &\mbox{on}&\partial\Omega
\end{array}
\right.
\end{align}
and similarly for $\beta$, with $\beta\geq u_0$, $0\geq \lambda c(x)\beta$ and reversed inequalities.
Therefore, theorem \ref{th2.1} gives us a solution $u_\lambda\in [\alpha,\beta]$, for all $\lambda<0$.

Observe that, since $\alpha\leq u_\lambda\leq\beta$, we can say that $\|u_\lambda\|_{C^{1,\alpha} (\overline{\Omega})}\leq C$ for all $\lambda\in [0,1]$, by the $C^{1,\alpha}$-estimates \cite{regularidade}. Thus, take a sequence $\lambda_k\leq 0$ with $\lambda_k \rightarrow 0$ as $k\rightarrow\infty$. Next, the compact inclusion $C^{1,\alpha}(\overline{\Omega})\subset E$ gives us some $u\in E$ such that $u_k\rightarrow u$ in $E$, up to a subsequence.
Hence we can define, for each $\varphi\in W^{2,p}_{\mathrm{loc}}(\Omega)$,
\begin{align}\label{def g e gk proofs}
g_k(x):=F(x,u_k,D\varphi,D^2\varphi)+h(x)+\lambda_kc(x)u_k ,\;
g(x):=F(x,u,D\varphi,D^2\varphi)+h(x) .
\end{align}
Then, $\|g_k-g\|_{L^p(\Omega)}\leq |\lambda_k|\|c\|_{L^p(\Omega)}\,\|u_0\|_\infty+\|d\|_{L^p(\Omega)}\,\omega (\|u_k-u\|_{L^\infty(\Omega)})\rightarrow 0$ as $k\rightarrow\infty$.
By proposition \ref{Lpquad}, we have that $u$ is an $L^p$-viscosity solution of $(P_0)$. From the uniqueness of the solution at $\lambda=0$, $u$ needs to be equal to $u_0$. Since the sequence of $\lambda$ converging to zero is arbitrary, we obtain $\|u_\lambda -u_0\|_E \rightarrow 0$ as $k\rightarrow\infty$.

\vspace{0.3cm}
Now we prove the existence of a continuum from $u_0$.
Fix an $\varepsilon>0$ and consider another pair of sub and supersolutions $\alpha :=u_0-\varepsilon$ and $\beta:=u_0+\varepsilon$. Analogously to \eqref{eq1 teo1}, we see that $\alpha,\beta$ are a pair of strong sub and supersolutions for $(P_0)$. Notice that they are not a pair for \ref{Plambda} with $\lambda<0$, since they do not have a sign. However, $\alpha<u_0<\beta$ in $\overline{\Omega}$, which implies that $\alpha\ll u_0\ll\beta$ in $\Omega\,$.
Since $u_0$ is the unique $L^p$-viscosity solution of the problem $(P_0)$, then $\alpha,\beta$ are strict in the sense of definition \ref{def2.1}. Then, theorem \ref{th2.1} $(i)$ gives us, for
$\mathcal{S}=\mathcal{O}\cap B_R$ defined there, that
\begin{align}\label{grau 0}
\mathrm{deg}(I-\mathcal{T}_0\,,\mathcal{S},\,0)=1.
\end{align}

Using again that $u_0$ is the unique $L^p$-viscosity solution of $(P_0)$, we have further that $\mathrm{ind}(I-\mathcal{T}_0\,,u_0)=1$. Then, by the well known degree theory results, see for example theorem 3.3 in \cite{BR}, there exists a continuum $\mathcal{C}\subset\Sigma$ such that both
$$
\mathcal{C}\cap ([0,+\infty)\times E) \;\;\;\mathrm{and}\;\;\;\mathcal{C}\cap ((-\infty,0]\times E)
$$
are unbounded in $\real^\pm \times E$. This proves item \textit{1.}

\vspace{0.3cm}
From now on, we suppose \ref{SC0}.

Let us prove point \textit{2.} in theorem \ref{th1.1,1.2,1.3}.
The continuum $\mathcal{C}\subset\Sigma$ is such that its projection on the $\lambda$-axis is either $\real$ (and we obtain $(ii)$ in theorem \ref{th1.1,1.2,1.3}) or it is $(-\infty,\bar{\lambda}]$, with $0<\bar{\lambda}<+\infty$. In the second case, since we know that the component $\mathcal{C}^+$ is unbounded in $R^+\times E$, its projection on the $E$ axis must be unbounded in $E$.

Under \ref{SC0}, by what we proved in the previous section, for any $0<\Lambda_1<\Lambda_2$ there is an $L^\infty$ a priori bound for the solutions of \ref{Plambda}, for all $\lambda\in [\Lambda_1,\Lambda_2]$. Then, by $C^{1,\alpha}$ global estimate \cite{regularidade}, we have also a $C^{1,\alpha}$ a priori bound for these solutions i.e the projection of $\mathcal{C}^+\cap ([\Lambda_1,\Lambda_2]\times E)$ on $E$ is bounded. So, $\mathcal{C}^+$ needs to be unbounded in $E$ when we approach $\lambda=0$ from the right.

Now, by proposition \ref{aprioribelowLinfty}, there is a lower $L^\infty$ bound for the solutions, for every $\lambda\leq\Lambda_2$. Therefore, $\mathcal{C}^+$ must emanate from plus infinity to the right of $\lambda =0$, with the positive part of its solutions blowing up to infinity in $C (\overline{\Omega})$. Thus, $(i)$ and $(ii)$ in \textit{2.} are proved.

\vspace{0.3cm}

Now we pass to the multiplicity results in item \textit{3.} of theorem \ref{th1.1,1.2,1.3}.

Observe that, up to taking a larger $R$ in \eqref{grau 0}, by $C^{1,\alpha}$-estimates we can suppose that
\begin{align}\label{u<R para lambda [0,1]}
\|u\|_{C^{1,\alpha}(\overline{\Omega})}<R, \textrm{ for all } L^p \textrm{-visc. solution } u \textrm{ of \ref{Plambda} in } [u_0-\varepsilon,u_0+\varepsilon], \; \lambda\in [0,1].
\end{align}

\begin{claim} \label{step 2 proof th1.3}
There exists a $\lambda_0>0$ such that $\,\mathrm{deg}(I-\mathcal{T}_\lambda\,,\mathcal{S},\,0)=1\,$, for all $ \lambda\in (0,\lambda_0)$.
\end{claim}

\begin{proof}
Let us prove the existence of a $\lambda_0>0$ such that, for all $\lambda\in (0,\lambda_0)$, \ref{Plambda} has no solution on $\partial\mathcal{S}$.
Suppose not, i.e. that for all $\lambda_0>0$, there exists a $\lambda\in (0,\lambda_0)$ such that $\mathcal{T}_\lambda$ has a fixed point on $\partial\mathcal{S}$.
Then for every $k\in\n$, there exists $\lambda_k\in \left(0,\frac{1}{k}\right)$ and $u_k\in \partial\mathcal{S}$ a solution of $(P_{\lambda_k})$. By \eqref{u<R para lambda [0,1]}, $u_k\not\in \partial B_R$ and so $u_k\in \partial\mathcal{O}$, for all $k\in\n$. Note that
$$
\partial\mathcal{O}=\{\,u\in C_0^1(\overline{\Omega})\,; \; \alpha\leq u\leq \beta \;\, \mathrm{in}\;\, \overline{\Omega}\;\;\mathrm{and}\,\; u \;\mathrm{``touches"}\;\alpha\;\mathrm{or}\;\beta\,\}
$$
where ``touches", as in \cite{BR}, has the following meaning.
\begin{defin}
Let $u,v\in C_0^1(\overline{\Omega})$. We say that $u$ ``touches" $v$ if $\exists\, x\in\Omega\,$ with $u(x)=v(x)$ or $\exists\, x\in\partial\Omega\,$ with $\partial_\nu u(x)= \partial_\nu v(x)$. In any case, $u(x)=v(x)$ at a point $x\in\,\overline{\Omega}$.
\end{defin}

If $u_k$ ``touches" $\alpha$, there exists a $x \in \overline{\Omega}\,$ such that $u_k(x)=u_0(x)-\varepsilon$, and since $u_k\geq u_0-\varepsilon$, then $\max_{\overline{\Omega}}\,(u_0-u_k)=\varepsilon$. If by other side $u_k$ ``touches" $\beta$, there exists a $x \in \overline{\Omega}\,$; $u_k(x)=u_0(x)+\varepsilon$, and since $u_k\leq u_0+\varepsilon$, then $\max_{\overline{\Omega}}\,(u_k-u_0)=\varepsilon$. Anyway,
\begin{align} \label{max|uk-u0|=epsilon}
\|u_k-u_0\|_\infty=\max_{\overline{\Omega}}\,|u_k-u_0|=\varepsilon\, , \;\;\textrm{ for all }\, k\in\n.
\end{align}

By \eqref{u<R para lambda [0,1]} and compact inclusion $C^{1,\alpha}(\overline{\Omega})\subset E$, there exists $u\in E$ such that $u_k\rightarrow u$ in $E$ as $k\rightarrow\infty$, up to a subsequence.
Hence, by the stability (proposition \ref{Lpquad}), $u$ is an $L^p$-viscosity solution of $(P_0)$. Indeed, we define $g$ and $g_k$ as in \eqref{def g e gk proofs}, for each $\varphi\in W^{2,p}_{loc}(\Omega)$, from which $\|g_k-g\|_{L^p(\Omega)}\leq \lambda_k\|c\|_{L^p(\Omega)}\,(\|u_0\|_\infty+\varepsilon)\rightarrow 0$ as $k\rightarrow\infty$. By uniqueness in $\lambda=0$, we have $u=u_0$, which contradicts \eqref{max|uk-u0|=epsilon} by taking the limit.

Therefore, the following degree is well defined and by the homotopy invariance and \eqref{grau 0} we obtain
$$
\mathrm{deg}(I-\mathcal{T}_\lambda\,,\mathcal{S},\,0)=\mathrm{deg}(I-\mathcal{T}_0\,,\mathcal{S},\,0)=1 , \, \textrm{ for all }\lambda\in (0,\lambda_0).
$$
\qedhere{\,\textit{Claim \ref{step 2 proof th1.3}.}}
\end{proof}

\begin{claim}\label{step 3 proof th1.3}
\ref{Plambda} has two solutions when $\lambda\in (0,\lambda_0/2]$.
\end{claim}

\begin{proof}
By claim \ref{step 2 proof th1.3}, the existence of a first solution $u_{\lambda,1}$ with $u_0-\varepsilon\ll u_{\lambda,1}\ll u_0+\varepsilon$ is already proved.
Set $\Lambda_2:=\lambda_0/2$.
Then, lemma \ref{lema Plambda,k has no solutions} implies that $(P_{\lambda,k})$ has no solutions for $k\geq 1$ and $\lambda\in (0,\Lambda_2]$.

Fix a $\lambda\in (0,\Lambda_2]$. With $h$ replaced by $h+k\widetilde{c}$  (see \eqref{ctilde e A for k=1}) we have, by theorem \ref{aprioriLinfty}, an $L^\infty$ a priori bound for solutions of $(P_{\lambda,k})$, for every $k\in [0,1]$. Precisely, we get an $L^\infty$ a priori bound for solutions of $(P_{\mu,k})$, for all $\mu\in [\lambda,\Lambda_2]$, depending on $\lambda$ and $\Lambda_2$. This provides, by the $C^{1,\alpha}$-estimates \cite{regularidade}, an a priori bound for solutions in $E$, namely
$$
\|u\|_{E} <R_0\, , \;\textrm{ for every } u\;L^p\textrm{-viscosity solution of } (P_{\lambda,k}), \textrm{ for all }\, k\in [0,1],
$$
and for some $R_0>R$ that depends, in addition to the coefficients of the equation, also on $\lambda$ and the $L^p$-norm of $\,\widetilde{c}$.
By the homotopy invariance of the degree in $k\in [0,1]$ and the fact that for $k=1$ the problem $(P_{\lambda,k})$ has no solution,
$$
\mathrm{deg}(I-\mathcal{T}_{\lambda}\,,B_{R_0}\,,\,0)=\mathrm{deg}(I-\mathcal{T}_{\lambda,0}\,,B_{R_0}\,,\,0)=\mathrm{deg}(I-\mathcal{T}_{\lambda,1}\,,B_{R_0}\,,\,0)=0
$$
where $\mathcal{T}_{\lambda,k}$ is the operator $\mathcal{T}_\lambda$ in which we replace $h$ by $h+k\widetilde{c}$ (of course $\mathcal{T}_{\lambda,k}$ keeps being completely continuous). But then, by the excision property of the degree,
\begin{align*}
\mathrm{deg}(I-\mathcal{T}_{\lambda}\,,B_{R_0}\setminus\mathcal{S},\,0)=\mathrm{deg}(I-\mathcal{T}_{\lambda}\,,B_{R_0}\,,\,0)-\mathrm{deg}(I-\mathcal{T}_\lambda\,,\mathcal{S},\,0)=-1
\end{align*}
which provides the second solution $u_{\lambda,2}\in B_{R_0}\setminus\mathcal{S}\,$ that we were looking for.

\qedhere{\,\textit{Claim \ref{step 3 proof th1.3}.}}
\end{proof}

Next, by claim \ref{step 3 proof th1.3}, the quantity
\begin{align*}
\bar{\lambda}:=\sup \{ \,\mu\, ; \;\forall\, \lambda\in (0,\mu), \; (P_\lambda )\;\, \mathrm{has\;at\;least\;two\;solutions}\} \in (0,+\infty]
\end{align*}
it is well defined and greater or equal than $\lambda_0/2$.

\begin{claim}\label{step 7 proof th.1.3}
$u_{\lambda,1}\rightarrow u_0$ in $E$ and $\max_{\overline{\Omega}}\, u_{\lambda,2}\rightarrow +\infty$ as $\lambda\rightarrow 0^+$.
\end{claim}

\begin{proof}
Let $(\lambda_k)_{k\in\n} \subset (0,\bar{\lambda})$ be a decreasing sequence with $\lambda_k\rightarrow 0$, so $\lambda_k \leq\lambda_0/2$ for $k\geq k_0$. Since $u_{\lambda_k,1}\in \mathcal{S}$, then $u_0-\varepsilon\leq u_{\lambda_k,1}\leq u_0+\varepsilon$ in $\overline{\Omega}$, therefore is bounded in $C^{1,\alpha}(\overline{\Omega})$ by \cite{regularidade}. Hence, exactly as in \eqref{def g e gk proofs}, we show by stability that $u_{\lambda_k,1}\rightarrow u$ in $E$, where $u$ is a solution of $(P_0)$. Therefore, $u=u_0$.

If, in turn, the respective sequence $u_{\lambda_k,2}$ were uniformly bounded from above, it would be unifomly bounded in $C(\overline{\Omega})$ using proposition  \ref{aprioribelowLinfty}, so bounded in $C^{1,\alpha}(\overline{\Omega})$ and the paragraph above would imply that $u_{\lambda_k,2}\rightarrow u_0$ in $E$. Since $u_0\in\mathcal{S}$ and $\mathcal{S}$ is open in $E$, then $u_{\lambda_k,2}$ should belong to $B_{r}^{E}(u_0)\subset\mathcal{S}$ for large $k$, for some $r>0$. But this contradicts the fact that $u_{\lambda_k,2}\notin \mathcal{S}$.
\qedhere{\,\textit{Claim \ref{step 7 proof th.1.3}.}}
\end{proof}

\begin{claim} \label{step 5 proof th1.3}
In case $\bar{\lambda}<+\infty$, the problem $(P_{\bar{\lambda}})$ has at least one solution.
\end{claim}

\begin{proof}
Let $\lambda_k\in (0,\bar{\lambda})$ be such that $\lambda_k\rightarrow \bar{\lambda}$ and let $u_k$ be a sequence of solutions for $(P_{\lambda_k})$.
Say $\lambda_k\in [\,\bar{\lambda}/2\,,\bar{\lambda}\,]$ for $k\geq k_0$. This provides an $L^\infty$ a priori bound for $u_k$, by theorem \ref{aprioriLinfty}, i.e. $\|u_k\|_\infty \leq C_1$, which implies that $\|u_k\|_{C^{1,\alpha}(\overline{\Omega})}\leq C_2$. Again, by compact inclusion and stability, we obtain $u_k\rightarrow u$ in $E$, where $u$ is a solution of $(P_{\bar{\lambda}})$.
Surely, for stability we need to consider, this time, $g(x):=F(x,u,D\varphi,D^2\varphi)+\bar{\lambda}c(x)u+h(x)$ and so $\|g_k-g\|_{L^p(\Omega)}\leq (\bar{\lambda}-\lambda_k)\,\|c\|_{L^p(\Omega)} \,C_2+\bar{\lambda}\,\|c\|_{L^p(\Omega)} \|u_k-u\|_\infty\rightarrow 0$ as $k\rightarrow\infty$.

\qedhere{\,\textit{Claim \ref{step 5 proof th1.3}.}}
\end{proof}

\vspace{0.3cm}
To finish the proof of theorem \ref{th1.1,1.2,1.3}, it remains to show the last statements in item~\textit{4.} concerning  ordering and uniqueness considerations, in which we assume \eqref{Hstrong}.
Notice that this automatically implies that solutions $u_{\lambda,1}$ and $u_{\lambda,2}$ are strong, as well as every $L^p$-viscosity solution of \ref{Plambda}.

Further, note that \eqref{Hstrong} provided existence of minimal and maximal solutions on the order interval $[\alpha,\beta]$ in theorem \ref{th2.1}. Such existence of minimal solution made it possible to find a minimal strong subsolution $\alpha_\lambda$ of \ref{Plambda} in lemma \ref{lemma 4.2}.
\vspace{0.3cm}

\begin{claim} \label{step 4 proof th1.3}
$u_{\lambda,1} \ll u_{\lambda,2}$, for all $\lambda\in (0,\bar{\lambda})$.
\end{claim}

\begin{proof}
Fix a $\lambda\in (0,\bar{\lambda})$ and consider the strict strong subsolution $\alpha=\alpha_\lambda$ given by lemma~\ref{lemma 4.2}.
Since in particular $\alpha\leq u$ for every (strong) solution of \ref{Plambda}, we can choose $u_{\lambda,1}$ as the minimal strong solution such that $u_{\lambda,1}\geq\alpha$ in $\Omega$. This choice implies that
\begin{align} \label{ulambda,1 leq ulambda,2 th1.3}
u_{\lambda,1}\lneqq u_{\lambda,2}\;\;\;\mathrm{in}\;\;\Omega.
\end{align}
Indeed, $u_{\lambda,1}\neq u_{\lambda,2}$ and, if there would exist $x_0\in\Omega$ such that $u_{\lambda,1}(x_0)>u_{\lambda,2}(x_0)$, by defining $u_\lambda:=\min_{\overline{\Omega}} \,\{ u_{\lambda,1}, u_{\lambda,2} \}$, as the minimum of strong supersolutions greater or equal than $\alpha$, so $u_\lambda\geq \alpha$ in $\Omega$. Therefore, theorem \ref{th2.1} would give us a solution $u$ of \ref{Plambda} such that $\alpha\leq u \leq u_\lambda \lneqq u_{\lambda,1}$, which contradicts the minimality of $u_{\lambda,1}$ and implies \eqref{ulambda,1 leq ulambda,2 th1.3}.

To finish the proof, define $v:=u_{\lambda,2}-u_{\lambda,1}\gneqq 0$ in $\Omega$ by \eqref{ulambda,1 leq ulambda,2 th1.3}. Then, since $u_{\lambda,1}$ and $u_{\lambda,2}$ are strong, $v$ satisfies, almost everywhere in $\Omega$,
\begin{align*}
-\mathcal{L}^-[v]\geq& -F[u_{\lambda,2}]+F[u_{\lambda,1}]\\=&\,\lambda c(x) v +\langle M(x)D v,D v \rangle+\langle M(x)D v,D u_{\lambda,1} \rangle +\langle M(x)D u_{\lambda,1},D v \rangle \\
\geq & -2\mu_2 |Du_{\lambda,1}|\, |Dv|.
\end{align*}
Hence, $v$ is a nonnegative strong solution of
$\mathcal{M}^-(D^2v)-\widetilde{b}(x)|Dv|\leq 0\,$ in $\Omega$, for $\widetilde{b}=b+2\mu_2 |Du_{\lambda,1}|\in L^p_+(\Omega)$.
Then SMP gives us that $v>0$ in $\Omega$, since $v\not\equiv 0$. Now, Hopf lemma and $v=0$ on $\partial\Omega$ imply $\partial_\nu v |_{\partial\Omega}>0$ and so $v\gg 0$ on $\Omega$.
\qedhere{\,\textit{Claim \ref{step 4 proof th1.3}.}}
\end{proof}

As far as uniqueness is concerned, from theorem 1(iii) in \cite{arma2010}, if the coercive problem for $\lambda\leq 0$ has a strong solution $u_\lambda$, it is the unique $L^p$-viscosity solution of \ref{Plambda}. So, under \eqref{Hstrong}, $u_\lambda$ is strong, then unique, in the $L^p$-viscosity sense, for all $\lambda<0$.
Observe that, in this case, we must have $\{\,(\lambda, u_\lambda), \; \lambda\leq 0\,\} \subset \mathcal{C}$. In another words, the projection of $\mathcal{C}$ on the $\lambda$-axis contains $(-\infty,0]$, as in theorems 1.1 and 1.2 in \cite{ACJT} for the Laplacian.

We finish the proof of \textit{4.} with the following claim.

\begin{claim}\label{uniqueness Pbarlambda th geral}
If $\bar{\lambda}<+\infty$ and $F$ is convex in $(r,p,X)$, the solution $u_{\bar{\lambda}}$ of $(P_{\bar{\lambda}})$, obtained in claim \ref{step 5 proof th1.3}, is unique.
\end{claim}

\begin{proof} Suppose, in order to obtain a contradiction, that there exist two different solutions $u_1$ and $u_2$ of $(P_{\bar{\lambda}})$, both strong by \eqref{Hstrong}. Consider $\beta=\beta_{\bar{\lambda}}:=\frac{1}{2} (u_1+u_2)$. Then, a.e. in $\Omega$,
\begin{align*}
-F\,[\beta]&\geq  -\{ F\,[u_1]+F\,[u_2]\}/2 \\
&= \,\bar{\lambda} c(x)\beta +h(x)+ \{\langle M(x) D u_1,D u_1\rangle +\langle M(x)D u_2,D u_2\rangle \}/2\\
&\gneqq\,\bar{\lambda} c(x)\beta +h(x)+\langle M(x) D\beta,D \beta\rangle
\end{align*}
using also the convexity of $\,\vec{p}\mapsto \langle M(x)\vec{p},\vec{p}\rangle$. Hence $\beta$ is a strong supersolution of $(P_{\bar{\lambda}})$ which is not a solution. Let us see that it is strict. Set $U:=\beta - u$, where $u\in E$ is an $L^p$-viscosity subsolution of $(P_{\bar{\lambda}})$ with $u\leq\beta$ in $\Omega$. Thus, $U$ is an $L^p$-viscosity solution of
\begin{align*}
-\mathcal{L}^-[U]&\gneqq \,\bar{\lambda} c(x) U -\langle M(x)D U,D U \rangle+\langle M(x)D U,D \beta \rangle +\langle M(x)D \beta,D U \rangle \\
&\geq \,-\mu_2|DU|^2 -2\mu_2 |D\beta|\, |DU|,
\end{align*}
and so, by lemma \ref{lemma2.3arma}, the function $w:=\frac{1}{m} (1-e^{-mU})$, where $m=\mu_2/\lambda_P$, is a nonnegative $L^p$-viscosity solution of
$\mathcal{M}^-(D^2w)-\widetilde{b}(x)|Dw| \lneqq 0\,$ in $\Omega$, with $\widetilde{b}=b+2\mu_2 |D\beta|\in L^p_+(\Omega)$.
Then SMP gives us that $w>0$ in $\Omega$. Now, Hopf and $w\geq 0$ on $\partial\Omega$ imply $\partial_\nu v |_{\partial\Omega}>0$ in the boundary points where $w=0$, and so $w\gg 0$ on $\Omega$. Consequently, $U\gg 0$ in $\Omega$ and $\beta$ is a strict strong supersolution of $(P_{\bar{\lambda}})$.

Consider also $\alpha=\alpha_{\bar{\lambda}}\,$ the strict strong subsolution of $(P_{\bar{\lambda}})$ given by proposition \ref{lemma 4.2} and look to the set $\bar{\mathcal{O}}=\{ \alpha\ll u\ll \beta \}=\{ \alpha_{\bar{\lambda}}\ll u\ll \beta_{\bar{\lambda}} \}$.
By the $C^{1,\alpha}$-estimates in \cite{regularidade},
\begin{align}\label{est th 1.3 Fconvex unicidade barlambda}
\|u\|_{C^{1,\alpha}(\overline{\Omega})}\leq C \;\textrm{ for all }\, u\in [\alpha,\beta], \;\; L^p \textrm{-visc. solution of } (P_\lambda), \textrm{ for all  } \lambda\in[\bar{\lambda},\bar{\lambda}+1].
\end{align}
for some $C>0$ that depends on the $L^\infty$-norm of $\alpha$. Then, by theorem \ref{th2.1}, we obtain $R>C$ such that $\mathrm{deg}(I-\mathcal{T}_{\bar{\lambda}}\,,\,\bar{\mathcal{S}},\,0)=1$, where $\bar{\mathcal{S}}=\bar{\mathcal{O}}\cap B_R$.

We claim that there exists $\varepsilon >0$ such that
\begin{align}\label{uniqueness barlambda eq degree epsilon}
\mathrm{deg}(I-\mathcal{T}_\lambda\,,\,\bar{\mathcal{S}},\,0)=1, \;\;\textrm{ for all } \lambda\in [\,\bar{\lambda},\bar{\lambda}+\varepsilon].
\end{align}
As in the proof of claim \ref{step 2 proof th1.3}, we will verify that there exists some $\varepsilon\in (0,1)$ such that there is no fixed points of $\mathcal{T}_\lambda$ on the boundary of $\bar{\mathcal{S}}$, for all $\lambda$ in the preceding interval.
Indeed, if this were not the case, there would exist a sequence $\lambda_k\rightarrow\bar{\lambda}$ with the respective solutions $u_k$ of $(P_{\lambda_k})$ belonging to $\partial\bar{\mathcal{S}}$.
Say $\lambda_k\in [\bar{\lambda},\bar{\lambda}+1]$ for $k\geq k_0$. Then, since $\alpha\leq u_k\leq \beta$ in $\Omega$, by \eqref{est th 1.3 Fconvex unicidade barlambda} we must have $u_k\in \partial\bar{\mathcal{O}}$ for $k\geq k_0$, which means that for each such $k$,
\begin{align}\label{th 1.3 contradiction unicidade em barlambda}
\max_{\overline{\Omega}} \,(\alpha-u_k)=0 \;\textrm{ or }\;\max_{\overline{\Omega}}\,(u_k-\beta)=0.
\end{align}
By \eqref{est th 1.3 Fconvex unicidade barlambda} and the compact inclusion $C^{1,\alpha}(\overline{\Omega}) \subset E$, $u_k\rightarrow u $ in $E$ for some $u\in E$, up to a subsequence.
This $u$ is an $L^p$-viscosity solution of $(P_{\bar{\lambda}})$ by the stability proposition \ref{Lpquad}; and $\alpha\leq u\leq \beta$ in $\Omega$ by taking the limit as $k\rightarrow +\infty$ in the corresponding inequalities for $u_k$.
Thus $\alpha\ll u \ll \beta$ in $\Omega$, since $\alpha$ and $\beta$ are strict. Passing to the limit into \eqref{th 1.3 contradiction unicidade em barlambda}, we obtain that $u$ ``touches" $\alpha$ or $\beta$, which contradicts the definition of $\alpha\ll u \ll \beta$.

Hence, obtaining \eqref{uniqueness barlambda eq degree epsilon} is just a question of applying homotopy invariance in $\lambda$ in the interval $[\bar{\lambda},\bar{\lambda}+\varepsilon]$.
Next, with \eqref{uniqueness barlambda eq degree epsilon} at hand, we repeat exactly the same argument done in claim \ref{step 3 proof th1.3} to obtain the existence of a second solution $u_{{\lambda},2}$ of \ref{Plambda}, for all $\lambda\in [\bar{\lambda},\bar{\lambda}+\varepsilon]$. But this, finally, contradicts the definition of $\bar{\lambda}$.
\qedhere{\,\textit{Claim \ref{uniqueness Pbarlambda th geral}.}}
\end{proof}

\subsection{Proof of Theorem \ref{th1.5}}

Suppose $u_0\leq 0$ with $cu_0\lneqq 0$ in $\Omega$ and \eqref{Hstrong}.

\begin{claim} \label{step 1 th1.5}
$u_0$ is a strict strong supersolution of \ref{Plambda}, for all $\lambda>0$.
\end{claim}

\begin{proof}
Since $\lambda c(x)u_0\lneqq 0$ in $\Omega$, $u_0$ is a strong supersolution of \ref{Plambda}  which is not a solution. To see that it is strict, we take $u\in E$ an $L^p$-viscosity subsolution of \ref{Plambda} such that $u\leq u_0$ in $\Omega$, and set $U:=u_0-u$. Then, since $u_0$ is strong, $U$ is an $L^p$-viscosity supersolution of
\begin{align*}
-\mathcal{L}^-[U]&\geq \lambda c(x) U -\langle M(x)D U,D U\rangle +\langle M(x)D u_0,D U\rangle+\langle M(x)D u_0,D U\rangle\\
&\geq - \mu_2\, |DU|^2-2\mu_2\, |Du_0|\,|DU|
\end{align*}
and so $\mathcal{M}^-(D^2w)-\widetilde{b}(x)|Dw|\leq 0$ in $\Omega$ in the $L^p$-viscosity sense, for $\widetilde{b}=b+2\mu_2 \,|Du_0|\in L^p_+(\Omega)$ and $w=\frac{1}{m}\{1-e^{-mU}\}$, $m=\mu_2/\lambda_P$, using lemma \ref{lemma2.3arma}.
Using SMP and the fact that $u_0$ is not a solution of \ref{Plambda}, we have $w>0$ in $\Omega$. Since $w\geq 0$ on $\partial\Omega$, at points belonging to $\partial\Omega$ such that $w>0$ we are done. If in turn $x_0\in\partial\Omega$ is such that $w(x_0)=0$, then $\partial_\nu U(x_0)>0$ by Hopf. Thus $w\gg 0$ and so $U\gg 0$ in $\Omega$.
\qedhere{\,\textit{Claim \ref{step 1 th1.5}.}}
\end{proof}

We now prove that for all $\lambda>0$, \ref{Plambda} has at least two solutions, $u_{\lambda,1}$ and $u_{\lambda,2}$, with $u_{\lambda,1} \ll u_0$ and $u_{\lambda,1}\ll u_{\lambda,2}$.

Fix a $\lambda>0$. From lemma \ref{lemma 4.2} and step 1, we get a pair of strong strict sub and supersolutions, $\alpha=\alpha_{\lambda}$ and $u_0$, which implies, by theorem \ref{th2.1}, the existence of a first solution $u_{\lambda,1}\in\mathcal{S}$, where $\mathcal{S}=\mathcal{S}_\lambda=\{ u\in C_0^1 (\overline{\Omega}); \; \alpha\ll u\ll u_0\;$ in $\Omega\}\cap B_R$ for some $R>0$.

\begin{rmk}
We already know, from theorem \ref{th1.1,1.2,1.3}, that at least two solutions exist. Here we redefine, as in \cite{CJ}, the set $\mathcal{S}$ in order to obtain more precise informations about $\Sigma$. Notice that, with this new definition, we automatically have $u_{\lambda,1}\ll u_0\leq 0$.
\end{rmk}

Fix a $\lambda>0$ and set $\Lambda_2:=2\lambda$. As in the proof of claim \ref{step 3 proof th1.3}, we observe that, by lemma~\ref{lema Plambda,k has no solutions}, $(P_{\lambda,k})$ has no solutions for $k\geq 1$. Moreover, for $h$ replaced by $h+k\widetilde{c}$, theorem~\ref{aprioriLinfty} gives us an $L^\infty$ a priori bound for solutions of $(P_{\lambda,k})$ for every $k\in [0,1]$, which depends on $\lambda$. This provides, by the $C^{1,\alpha}$ global estimates \cite{regularidade}, an a priori bound for solutions in $E$, i.e. $\|u\|_{E} <R_0$ for every $u$ solution of  $(P_{\lambda,k})$, for all $k\in [0,1]$, where $R_0>R$ also depends on $\lambda$.
By the homotopy invariance of the degree,
$$
\mathrm{deg}(I-\mathcal{T}_{\lambda}\,,B_{R_0}\,,\,0)=\mathrm{deg}(I-\mathcal{T}_{\lambda,0}\,,B_{R_0}\,,\,0)=\mathrm{deg}(I-\mathcal{T}_{\lambda,1}\,,B_{R_0}\,,\,0)=0.
$$
Therefore, by the excision property of the degree
\begin{align*}
\mathrm{deg}(I-\mathcal{T}_{\lambda}\,,B_{R_0}\setminus\mathcal{S},\,0)=\mathrm{deg}(I-\mathcal{T}_{\lambda}\,,B_{R_0}\,,\,0)-\mathrm{deg}(I-\mathcal{T}_\lambda\,,\mathcal{S},\,0)=-1
\end{align*}
and the existence of a second solution $u_{\lambda,2}\in B_{R_0}\setminus\mathcal{S}\,$ is derived.

Since the argument above can be done for any $\lambda>0$, we obtain the existence of at least two solutions for every positive $\lambda$. Exactly the same reasoning from claims \ref{step 4 proof th1.3} and  \ref{step 7 proof th.1.3} applies to check that  $u_{\lambda,1} \ll u_{\lambda,2}$ in $\Omega$ and to get their behavior when $\lambda \rightarrow 0^+$, respectively, since $\bar{\lambda}$ is the same from theorem \ref{th1.1,1.2,1.3}. Of course, here $\bar{\lambda}=+\infty$.

\begin{claim} \label{step 6 proof th1.5}
For $\lambda_1<\lambda_2\,$, we have $u_{\lambda_2,1}\ll u_{\lambda_1,1}$ in $\Omega$.
\end{claim}

\begin{proof}
For fixed $\lambda_1<\lambda_2\,$ note that $\lambda_1 \,c(x)\,u_{\lambda_1,1} \gneqq \lambda_2 \,c(x)\,u_{\lambda_1,1}\,$ since $u_{\lambda,1}< 0$. Then, $u_{\lambda_1,1}$ is a strong supersolution of $(P_{\lambda_2})$ which is not a solution and, in particular, $u_{\lambda_1,1}\neq u_{\lambda_2,1}$.

We first infer that $u_{\lambda_2,1}\lneqq u_{\lambda_1,1}\,$ in $\Omega$.
In fact, similarly to the argument in the proof of claim \ref{step 4 proof th1.3}, recall that $\alpha=\alpha_{\lambda_2}$, given by lemma \ref{lemma 4.2}, is such that $\alpha\leq u$ for every strong supersolution of $(P_{\lambda_2})$, and in particular $\alpha \leq u_{\lambda_1,1}$. Remember also that $u_{\lambda_2,1}$ is the minimal strong solution such that $u_{\lambda_2,1}\geq\alpha$ in $\Omega$.
Now, if there was a $x_0\in\Omega$ such that $u_{\lambda_2,1}(x_0)>u_{\lambda_1,1}(x_0)$, by defining $\beta:=\min_{\overline{\Omega}} \,\{ u_{\lambda_1,1}, u_{\lambda_2,1} \}$, as the minimum of strong supersolutions of $(P_{\lambda_2})$ not less than $\alpha$, we have $\alpha\leq\beta$ in $\Omega$.
Thus, theorem \ref{th2.1} provides a solution $u$ of $(P_{\lambda_2})$ such that $\alpha\leq u \leq \beta \lneqq u_{\lambda_2,1}$ in $\Omega$, which contradicts the minimality of $u_{\lambda_2,1}$.

Proceeding as usual, $v:=u_{\lambda_1,1}-u_{\lambda_2,1}$ becomes a nonnegative strong supersolution of
$\mathcal{M}^-(D^2v)-\widetilde{b}(x)|Dv|\leq 0$ in $\Omega$,
then SMP gives us that $v>0$ in $\Omega$, since $v\not\equiv 0$; Hopf and $v=0$ on $\partial\Omega\,$ imply that $\partial_\nu v |_{\partial\Omega}>0$, so $v\gg 0$ on $\Omega$.
\qedhere{\,\textit{Claim \ref{step 6 proof th1.5}.}}
\end{proof}

\begin{rmk} \label{step 1 prop4.3}
Notice that $u\ll u_0$ in $\Omega$, for every nonpositive $L^p$-viscosity subsolution $u$ of \ref{Plambda} in $E$.
Indeed, since $\lambda c(x)u\leq 0$ in $\Omega$, $u$ is also an $L^p$-viscosity subsolution of $(P_0)$. By remark \ref{obs.c.leq0}, $u\leq u_0$, since $u_0$ is strong. Now, by claim \ref{step 1 th1.5} and definition of strict supersolution, we get $u\ll u_0$ in $\Omega$.

In particular, $u\equiv 0$ is never a solution of \ref{Plambda}, for any $\lambda>0$ $($recall that $cu_0\not\equiv 0)$.
\end{rmk}

\begin{claim}\label{prop4.3}
In addition to the hypothesis of theorem \ref{th1.5}, suppose that $F$ is convex in $(r,\vec{p},X)$. Then, \ref{Plambda} has at most one nonpositive solution. In particular, $\max_{\overline{\Omega}}\,u_{\lambda,2}>0$.
\end{claim}

\begin{proof}
Suppose, in order to obtain a contradiction, that there exist two different nonpositive solutions $u_1$ and $u_2$, strong by \eqref{Hstrong}. By remark \ref{step 1 prop4.3} we know that $u_1\ll u_0$ and $u_2\ll u_0$ in $\Omega$.
We can assume that they are ordered, in the sense that $u_1\lneqq u_2$. Indeed, observe that $\max \{u_1,u_2\} \leq u_0$, then by theorem \ref{th2.1} we obtain a solution $u_3$ of \ref{Plambda} satisfying $u_2\leq \max \{u_1,u_2\}\leq u_3 \leq u_0\leq 0$ in $\Omega$. Thus, if the solutions $u_1$ and $u_2$ do not satisfy  $u_1\lneqq u_2$, then there is a point $x_0\in\Omega$ with $u_1 (x_0)>u_2(x_0)$, which implies that $u_2\not\equiv \max \{u_1,u_2\}$ and so  $u_2\lneqq u_3$; in this case we just replace $u_1,u_2$ by $u_2,u_3$ respectively.

Since $u_2\ll0$ (from $u_2\ll u_0$ and $\partial_\nu u_0\leq 0$), the quantity $\tau:=\inf\{t>0; \, (1+t)u_2\leq u_1\textrm{ in }\Omega\}$ is well defined and finite. Further, $\tau>0$, since $u_2-u_1\gneqq 0$, so this infimum is attained. Then, by setting $w:=\frac{1}{\tau}\{(1+\tau)u_2-u_1\}$, we have that $w\leq 0$ satisfies $u_2=\frac{\tau}{1+\tau}w+\frac{1}{1+\tau}u_1$ and it is a strong subsolution of
\begin{align*}
F[w]&\geq \frac{1+\tau}{\tau} F[u_2]-\frac{1}{\tau}F[u_1]=
\frac{1+\tau}{\tau} \{\,\lambda c(x)u_2 +\langle M(x)D u_2,D u_2\rangle+h(x) \,\}\\ &-\frac{1}{\tau} \{\, \lambda c(x)u_1 +\langle M(x)D u_1,D u_1\rangle +h(x)\,\}\\
&\geq \lambda c(x)w +\langle M(x)D w,D w\rangle +h(x) \quad \textrm{ in } \Omega
\end{align*}
since $F$ is convex in $(r,\vec{p},X)$ and $\vec{p}\mapsto\langle M(x)\vec{p},\vec{p}\rangle$ is convex in $\vec{p}$. Now, by remark \ref{step 1 prop4.3} we have $w\ll u_0$ in $\Omega$, i.e. $w=w_\tau<0$ in $\Omega$. Then, there exists a little bit smaller $t\in (0,\tau)$ such that $w_t<0$ in $\Omega$ (see the argument in \cite{regularidade}, for example, by taking a compact set with small measure containing the boundary).
Therefore, this last contradicts the definition of $\tau$ as a minimum.
\qedhere{\,\textit{Claim \ref{prop4.3}.}}
\end{proof}

\subsection{Proof of Theorem \ref{th1.4}}

Suppose for the time being just $u_0\geq 0$ with $cu_0\gneqq 0$ in $\Omega$ and \ref{SC0}.

\begin{claim} \label{step 1 proof th1.4}
We have $u\gg u_0$, for every nonnegative $L^p$-viscosity supersolution $u\in E$ of \ref{Plambda}, for all $\lambda>0$.
\end{claim}

\begin{proof}
Notice that $\lambda c(x)u\geq 0$ in $\Omega$ implies that $u$ is an $L^p$-viscosity supersolution of $(P_0)$. Since $u_0$ is strong, by remark \ref{obs.c.leq0}, $u\geq u_0$ in $\Omega$. But $u_0$ is not a solution of \ref{Plambda} for $\lambda>0$ since $cu_0\not\equiv 0$, which means that $u_0\not\equiv u$.

Set $v:=u-u_0$ in $\Omega$. Then, using $M\geq 0$, we see that $v$ is a nonnegative $L^p$-viscosity supersolution of
$$
\mathcal{M}^-(D^2v)-\widetilde{b}(x)|Dv|\leq 0 \quad\mathrm{in}\;\;\Omega
$$
with $\widetilde{b}:=b+2\mu_2\,|Du_0|$, as usual.  By SMP, $v>0$ in $\Omega$. If $v>0$ on $\partial\Omega$ it is done; if by other side there exists $x_0\in \partial\Omega$ with $v(x_0)=0$, we apply Hopf lemma to obtain $\partial_\nu v(x_0)>0$. Therefore, $v\gg0$ in $\Omega$.
\qedhere{\,\textit{Claim \ref{step 1 proof th1.4}.}}
\end{proof}

\begin{claim} \label{no nonnegative sol of Plambda for lambda large}
\ref{Plambda} has no nonnegative $L^p$-viscosity solutions for $\lambda$ large.
\end{claim}

\begin{proof}
Let $\lambda\geq \widetilde{\lambda}_1$, where $\widetilde{\lambda}_1=\widetilde{\lambda}_1^+\, ( \widetilde{\mathcal{L}}^-(c),\Omega )>0$ is the principal weighted eigenvalue of
\begin{align*}
\widetilde{\mathcal{L}}^-\,[v]:=\mathcal{M}^-(D^2v)+\widetilde{b}(x)|Dv|\,,\;\;\widetilde{b}(x):=b(x)+2\mu_2\, |D u_0|\in L^p_+(\Omega) ,
\end{align*}
associated to $\widetilde{\varphi}_1=\widetilde{\varphi}_1^+\,( \widetilde{\mathcal{L}}^-(c),\Omega )\in W^{2,p}(\Omega)$, from proposition \ref{exist eig for L-c geq 0}, i.e.
\begin{align} \label{def phi1tilde+ of L-(c)}
\left\{
\begin{array}{rclcc}
(\widetilde{\mathcal{L}}^-+\widetilde{\lambda}_1c) [\widetilde{\varphi}_1] &=& 0 &\mbox{in}& \Omega \\
\widetilde{\varphi}_1 &>&0  &\mbox{in} &\Omega \\
\widetilde{\varphi}_1&=&0  &\mbox{on} &\partial\Omega
\end{array}
\right.
\end{align}

Suppose, in order to obtain a contradiction, that there exists a nonnegative $L^p$-viscosity solution $u$ of \ref{Plambda} and set $v:=u-u_0$ in $\Omega$. By claim \ref{step 1 proof th1.4}, $v\gg 0$ in $\Omega$.

Since $u_0$ is strong, we can use it as a test function into the definition of $L^p$-viscosity supersolution of $u$, together with \ref{SC0} and $\mu_2I \geq M\geq 0$, to obtain that
\begin{align*}
-\mathcal{L}^-[v]&\geq \lambda c(x)v+\lambda c(x) u_0 +\langle M(x)D v, D v\rangle +\langle M(x)D v, D u_0\rangle +\langle M(x)D u_0, D v\rangle \\
&\gneqq \widetilde{\lambda}_1 c(x) v-2\mu_2 |D u_0|\, |D v|
\end{align*}
since $c(x)u_0\gneqq0$, and so $v$ satisfies
\begin{align}\label{eq v L-(c) to absurd th1.4}
\left\{
\begin{array}{rclcc}
(\widetilde{\mathcal{L}}^-+\widetilde{\lambda}_1c)\, [v] &\lneqq & 0 &\;\mbox{in}&\Omega \\
v&>& 0 &\;\mbox{in} & \Omega
\end{array}
\right.
\end{align}
in the $L^p$-viscosity sense. As the proof of lemma \ref{lema Plambda,k has no solutions}, applying proposition \ref{th4.1 QB} to \eqref{eq v L-(c) to absurd th1.4} and \eqref{def phi1tilde+ of L-(c)}, we get $v=t\widetilde{\varphi}_1$ for $t>0$. But this contradicts the first line in \eqref{eq v L-(c) to absurd th1.4}, since $(\widetilde{\mathcal{L}}^-+\widetilde{\lambda}_1c) \,[t\widetilde{\varphi}_1]=0$ in $\Omega$.
\qedhere{\,\textit{Claim \ref{no nonnegative sol of Plambda for lambda large}.}}
\end{proof}

Define
$$
\bar{\lambda}:=\sup \,\{\, \lambda \,; \; (P_\lambda) \;\,\textrm{has an } L^p\textrm{-viscosity solution }u_\lambda\geq 0 \;\,\mathrm{in}\;\Omega \,\}
$$
which is finite, by claim \ref{no nonnegative sol of Plambda for lambda large}. Of course it is well defined and nonnegative, since $u_0\geq0$. Also, by the definition of $\bar{\lambda}$, \ref{Plambda} has no nonnegative solutions for $\lambda>\bar{\lambda}$.

It is a subtle but important detail that $\bar{\lambda}$ is a positive number. In fact, by the existence of the continuum, theorem \ref{th1.1,1.2,1.3}, we know that, for $\lambda$ small, there exists a solution $u_{\lambda}$ of $(P_{\lambda})$ such that $u_{\lambda}\not\equiv 0$, since $\|u_0\|_{E}>0$. But why can we infer that $u_{\lambda}\geq 0$ for small $\lambda$ positive? This is the subject of the next claim.

Consider hypothesis \eqref{Hstrong} from now. Then, $L^p$-viscosity solutions of \ref{Plambda} are strong.

\begin{claim}
$\bar{\lambda}>0$.
\end{claim}

\begin{proof}\label{barlambda positivo th1.4}
Let $\lambda\in (0,\Lambda_0)$, where $\Lambda_0$ is such that there exists a nontrivial solution of \ref{Plambda} in this interval, as indicated above. Suppose $\Lambda_0\leq \min \{1,1/C_0\}$, where $C_0$ is a lower bound for the solutions of \ref{Plambda} such that $u\geq -C_0$, for all $\lambda\in [0,1]$. We are supposing here $C_0>0$, otherwise every solution of \ref{Plambda} would be nonnegative for $\lambda\leq 1$.

Suppose firstly that $h\geq c$. In this case every nontrivial solution $u$ of \ref{Plambda} satisfies
\begin{align*}
\left\{
\begin{array}{rclcc}
-\mathcal{L}^-[u]&\geq &-F[u]\geq c(x)(1-\Lambda_0 \,C_0)+\langle M(x)D u,D u\rangle \geq 0 &\mbox{ in } & \Omega  \\
u&=&0 & \mbox{ on } &\partial\Omega
\end{array}
\right. .
\end{align*}
Then $u\geq 0$ in $\Omega$ by ABP, for all $\lambda\in (0,\Lambda_0)$.

In the general case, let $u_c$ a nontrivial nonnegative solution of the problem
\begin{align*}
\left\{
\begin{array}{rclcc}
-F[u_c] &=& \widetilde{h}(x)+\lambda c(x) u_c+\langle M(x)D u_c,D u_c \rangle \geq 0 &\mbox{ in } & \Omega \\
u_c &=& 0 &\mbox{ on } & \partial\Omega
\end{array}
\right.
\end{align*}
for $\widetilde{h}:=\max \{h,c\} \geq c$. Notice that $u_c$ is strong under \eqref{Hstrong}, since $|\widetilde{h}|\leq |h|+c$. Moreover, since $\widetilde{h}\geq h$, $u_c$ is a supersolution of \ref{Plambda}, then $u_c\geq u_0$ by claim \ref{step 1 proof th1.4}.

Further, since $\lambda c(x)u_0\geq 0$, $u_0$ is a strong subsolution of \ref{Plambda}.
Thus, applying theorem \ref{th2.1}, we obtain an $L^p$-viscosity solution $u$ of \ref{Plambda} with $u_0\leq u\leq u_c$ in $\Omega$. In particular this solution is nonnegative and nontrivial.
\qedhere{\,\textit{Claim \ref{barlambda positivo th1.4}. }}
\end{proof}

\begin{rmk}
Another way to prove claim \ref{barlambda positivo th1.4} is by considerations on first eigenvalues, which provide an  estimate on the smallness of $\lambda$.

Set $v:=u_0-u$, with $u$ an $L^p$-viscosity solution of \ref{Plambda} for $0<\lambda <(C_A\,\|c\|_{L^p(\Omega)})^{-1}$, where $C_A$ is the constant from ABP for $\mu=0$. Of course,  negativity of $v$  implies that $u\geq u_0\geq 0$ in $\Omega$. So, in order to obtain a contradiction, suppose that $\sup_{\Omega} v>0$.

Notice that $v$ is an $L^p$-viscosity solution of $\widetilde{\mathcal{L}}^+ [v]\geq -\lambda c(x)v^+$ in $\Omega$, with $v=0$ on $\partial\Omega$.
Thus, as in the proof of proposition 3.4 in \cite{arma2010}, we use ABP to obtain that $ \sup_{\Omega} v \leq \lambda\, C_A  \|c\|_{L^p (\Omega)}\sup_{\Omega} v^+ $ which, by the choice of $\lambda$, is a contradiction.
\end{rmk}

\begin{claim}\label{step 4 proof th1.4}
For each $\lambda\in (0,\bar{\lambda})$, \ref{Plambda} has a well ordered strict pair of strong sub and supersolutions, namely $u_0 \ll \beta_\lambda$ in $\Omega$.
\end{claim}

\begin{proof} Let $\lambda\in (0,\bar{\lambda})$. As the strict subsolution we just consider $u_0$ again, which is strong.
Note that $u_0$ is strict, since for any supersolution $u\in E$ such that $u\geq u_0$ in $\Omega$, we have $u\gg u_0$, by repeating the  final paragraph in the proof of claim \ref{step 1 proof th1.4}.

Note that, from the definition of $\bar{\lambda}$, there exists $\mu\in (\lambda,\bar{\lambda})$ and a nonnegative solution $u_{\mu}$ of $(P_\mu)$. By claim \ref{step 1 proof th1.4}, $u_\mu\gg u_0$ in $\Omega$.
On the other hand, since
$$c(x)(\,\mu-\lambda) \,u_\mu \geq c(x)(\,\mu-\lambda) \,u_0 \gneqq 0,$$
we have $\mu c(x)u_\mu \gneqq \lambda c(x)u_\mu$ and so $u_\mu$ is a supersolution of \ref{Plambda} which is not a solution.
In addition, $u_\mu$ is strict because if $u\in E$ is an $L^p$-viscosity subsolution of \ref{Plambda} with $u\leq u_\mu$ in $\Omega$, by defining $v=u_\mu -u$ and arguing as usual when we have a strong supersolution, $v$ becomes an $L^p$-viscosity supersolution of
$$
\mathcal{M}^-(D^2v)-\widetilde{b}(x)|Dv|-\mu_2\,|D v|^2\lneqq 0\;\;\textrm{ in }\Omega
$$
with $\widetilde{b}:=b+2\mu_2\,|Du_\mu|\in L^p_+(\Omega)$. Thus, $w=\frac{1}{m} \{1-e^{-mv}\}$, for $m=\mu_2/\lambda_P$, is a  $L^p$-viscosity supersolution of $\mathcal{M}^-(D^2w)-\widetilde{b}(x)|Dw|\lneqq 0$ in $\Omega$ by lemma \ref{lemma2.3arma}.
Then, SMP gives us $w>0$ in $\Omega$ and so $w\gg 0$ in $\Omega$, by applying Hopf at the boundary bounds where $w=0$. Consequently, $v \gg 0$ in $\Omega$.

Therefore, we can define $\beta_\lambda:=u_{\mu}$ for $\mu=\mu (\lambda)$, for all $\lambda\in (0,\bar{\lambda})$.

\qedhere{\,\textit{Claim \ref{step 4 proof th1.4}. }}
\end{proof}

Hence, by claim \ref{step 4 proof th1.4} and theorem \ref{th2.1}, there exists a solution $u_{\lambda,1}$ of \ref{Plambda} with $u_0\leq u_{\lambda,1} \leq \beta_\lambda\,$ in $\Omega$ and $\mathrm{deg} (I-\mathcal{T}_\lambda, \mathcal{S},0)=1$, for all $\lambda>0$.

Next, we work a little bit more to construct the second solution $u_{\lambda,2}$ that also satisfies $u_{\lambda,2}\gg u_0$ but is not on $\mathcal{S}$ (as in \cite{CJ}). For this, fix a $\lambda\in (0,\bar{\lambda})$ and consider the open subset of $E$ defined by $$\mathcal{D}=\{ u\in C_0^1(\overline{\Omega}); \; u_0\ll u  \; \mathrm{in}\;\Omega \},$$ which contains the set $\mathcal{S}$ from theorem \ref{th2.1}, since $\mathcal{S}=\{ u\in B_R; \, u\ll \beta_\lambda\; \mathrm{in}\;\Omega\}\cap \mathcal{D}$.

Analogously to the proof of claim \ref{step 3 proof th1.3}, we obtain an  a priori $L^\infty$-bound for the solutions of $(P_{\lambda,k})$ which depends on $\lambda$ but not on $k\in [0,1]$, related to $\Lambda_2:=\bar{\lambda}-\delta$ for some small $\delta >0$.
This provides a $R_0=R_0 \,(\lambda,\delta)>R$ which bounds the $E$ norm of the solutions, by the $C^{1,\alpha}$-estimates. Then, by the homotopy invariance of the degree in $k$ and the fact that there is no solution for $k=1$, we have $\mathrm{deg} (I-\mathcal{T}_\lambda, B_{R_0}\cap\mathcal{D}\,,0)=0$. Therefore, by excision, $\mathrm{deg} (I-\mathcal{T}_\lambda, (B_{R_0}\cap\mathcal{D})\setminus\mathcal{S},0)=-1$, which provides a second solution $u_{\lambda,2}\in\mathcal{D}\setminus\mathcal{S}$, i.e. a solution that satisfies, by construction, $u_{\lambda,2}\gg u_0$ in $\Omega$, for $\lambda\in (0,\bar{\lambda}-\delta]$, for every $\delta >0$. In particular, this second solution is also nonnegative and nontrivial for all $\lambda < \bar{\lambda}$.

\vspace{0.1cm}
Under \eqref{Hstrong}, theorem \ref{th2.1} (ii) allows us to choose $u_{\lambda,1}$ as the minimal strong solution between $u_0$ and $\beta_\lambda$.
As the proof of claim \ref{step 4 proof th1.3}, this implies \eqref{ulambda,1 leq ulambda,2 th1.3}.  Indeed, $u_{\lambda,1}\neq u_{\lambda,2}$ and, if existed $x_0\in\Omega$ such that $u_{\lambda,1}(x_0)>u_{\lambda,2}(x_0)$, by defining $\,\tilde{\beta}=\tilde{\beta}_\lambda:=\min_{\overline{\Omega}} \,\{ u_{\lambda,1}, u_{\lambda,2},\beta_\lambda\}$, as the minimum of strong supersolutions greater or equal than $u_0\,$, we have $\tilde{\beta}\geq u_0$ in $\Omega$. Also, $\tilde{\beta}\leq \beta_\lambda$ in $\Omega$.
By theorem \ref{th2.1} there exists a solution $u$ of \ref{Plambda} such that $u_0\leq u \leq \tilde{\beta} \lneqq u_{\lambda,1}$, which contradicts the minimality of $u_{\lambda,1}$, since $u$ is a solution which belongs to the order interval $[u_0,\beta_\lambda]$.

Therefore, defining $v=u_{\lambda,2}-u_{\lambda,1}\gneqq 0$ in $\Omega$, we see that $v$ is a nonnegative strong solution of
$\mathcal{M}^-(D^2v)-\widetilde{b}(x)|Dv|\leq 0\,$ in $\Omega$.
Thus, since $v\not\equiv 0$, SMP yields $v>0$ in $\Omega$ and Hopf concludes that $v\gg 0$ in $\Omega$, i.e.
\begin{center}
$u_{\lambda,1}\ll u_{\lambda,2}\,$ in $\,\Omega$, for all $\,\lambda\in (0,\bar{\lambda})$.
\end{center}

\begin{claim} \label{step 6 proof th1.4}
For $\lambda_1<\lambda_2\,$, we have $u_{\lambda_1,1}\ll u_{\lambda_2,1}$ in $\Omega$.
\end{claim}

\begin{proof} The proof is similar to the proof of claim \ref{step 4 proof th1.4}, but a little bit simpler since both $u_{\lambda_1,1}$ and $u_{\lambda_2,2}$ are strong. However, we repeat it here to avoid confusions about notation.
For fixed $\lambda_1<\lambda_2\,$, we have $\lambda_2 \,c(x)\,u_{\lambda_2,1} \gneqq \lambda_1 \,c(x)\,u_{\lambda_2,1}\,$ since
$$(\lambda_2-\lambda_1)c(x)u_{\lambda_2,1}\geq (\lambda_2-\lambda_1)c(x)u_0 \gneqq 0.$$
Then, $u_{\lambda_2,1}$ is a supersolution of $(P_{\lambda_1})$ that is not a solution. In particular, $u_{\lambda_1,1}\neq u_{\lambda_2,1}$.

Next $u_{\lambda_1,1}\lneqq u_{\lambda_2,1}\,$ in $\Omega$. In fact, if there was a point $x_0\in\Omega$ such that $u_{\lambda_1,1}(x_0)>u_{\lambda_2,1}(x_0)$, by defining $\tilde{\beta}=\min_{\overline{\Omega}} \,\{ u_{\lambda_1,1}, u_{\lambda_2,1},\beta_{{\lambda}_1} \}$, the minimum of strong supersolutions of $(P_{\lambda_1})$ larger than $u_0$, we have that $u_0\leq\tilde{\beta}$ in $\Omega$. By theorem \ref{th2.1}, there exists a solution $u$ of $(P_{\lambda_2})$ such that $u_0\leq u \leq \tilde{\beta} \lneqq u_{\lambda_1,1}$, which contradicts the minimality of $u_{\lambda_1,1}$, since $u$ is a solution that belongs to the order interval $[u_0,\beta_{{\lambda}_1}]$.

Hence $v:=u_{\lambda_2,1}-u_{\lambda_1,1}\gneqq 0$ is a strong supersolution of
$\mathcal{M}^-(D^2v)-\widetilde{b}(x)|Dv|\leq 0\,$ in $\Omega$, for $\widetilde{b}=b+2\mu_2\,|Du_{\lambda_1,1}|$.
Then SMP yields $v>0$ in $\Omega$, since $v\not\equiv 0$.  Hopf and $v=0$ on $\partial\Omega$ give us $\partial_\nu v |_{\partial\Omega}>0$, so $v\gg 0$ in $\Omega$.
\qedhere{\,\textit{Claim \ref{step 6 proof th1.5}.}}
\end{proof}

The existence proof of at least one solution for $(P_{\bar{\lambda}})$ follows exactly the same lines as the proof of claim \ref{step 5 proof th1.3} since, in there, we only used the fact that there exists one sequence of solutions corresponding to a maximizing sequence of $\lambda$'s converging to the supremum $\bar{\lambda}$. Furthermore, uniqueness is true if $F$ is convex in $(r,p,X)$, by following the proof of claim \ref{uniqueness Pbarlambda th geral}.
Finally, the behavior of the solutions is the same as in claim \ref{step 7 proof th.1.3} and this finishes the proof of theorem \ref{th1.4}.


\begin{rmk}\label{particularh}
Particular cases of theorems \ref{th1.5} and \ref{th1.4} are $h\lneqq 0$ and $h\gneqq 0$, respectively.
Indeed,  if $h\gneqq 0$ and \ref{H0} holds, then $u_0$ is a strong supersolution of
$$
-\mathcal{L}^-[u_0]\geq -F[u_0]=h(x)+\langle M(x)D u_0,D u_0 \rangle \gneqq 0 \quad \textrm{in}\;\;\Omega
$$
with $u_0=0$ on $\partial\Omega$. Then SPM gives us $u_0>0$ in $\Omega$ and so $cu_0\not\equiv 0$. Furthermore, by Hopf, $u_0\gg 0$ in $\Omega$.

On the other hand, if $h\lneqq 0$ and \ref{H0} holds, then $u_0$ is a strong subsolution of
$$
-\mathcal{L}^+[u_0]\leq -F[u_0]\lneqq \langle M(x)D u_0,D u_0 \rangle  \quad \textrm{in}\;\;\Omega
$$
and so $v_0:=\frac{1}{m}(e^{mu_0}-1)$, for $m=\frac{\mu_2}{\lambda_P}$, is a strong subsolution of
$\mathcal{L}^+[v_0]\gneqq 0$ in $\Omega$ by lemma \ref{lemma2.3arma}, with $v_0=0$ on $\partial\Omega$. Again by SMP we get $v_0<0$ in $\Omega$ $($then $v_0\ll0$ in $\Omega$ by Hopf $)$ and so does $u_0$ $($with $u_0\ll 0)$, from where $cu_0\not\equiv 0$.

Notice that in the case $h\equiv 0$, we have that $u\equiv 0$ is a strong solution of \ref{Plambda}, for all $\lambda\in\real$. By theorem 1(iii) in \cite{arma2010}, this is the unique $L^p$-viscosity solution for all $\lambda\leq 0$. By theorem \ref{th1.1,1.2,1.3} we obtain the existence of $\bar{\lambda}$ such that \ref{Plambda} has at least one more nontrivial solution, for all $\lambda\in (0,\bar{\lambda})$.
\end{rmk}


\section{Appendix}

\subsection{Proof of Lemma \ref{vazquez particular} }\label{appendix1}

For maximal generality we  show that the result is valid for the most general  notion of $C$-viscosity solutions (see \cite{user}, \cite{Koike}), which can be discontinuous functions.
We can suppose $u\in LSC(\Omega)$, just replacing $u$ by its lower semi-continuous envelope $u_*$, defined as $u_*(x)=\varliminf_{y\rightarrow x}\, u(y)\in LSC(\Omega)$. In this case, we say that $u$ is a $C$-viscosity supersolution if $u_*$ is.

By contradiction, let $u$ be a nonnegative $C$-viscosity supersolution of
$$
\mathcal{L}_1^- [u]\leq f(u) \quad \mathrm{in}\;\;\Omega
$$
with both $\Omega_0:=\{x\in \Omega; \;u(x)=0\}$ and $\Omega^+:=\{x\in\Omega;\;u(x)>0\}$ nonempty sets. Notice that $\Omega^+$ is open, since $u\in LSC(\Omega)$. As in the usual proof of SMP (see, for instance, \cite{BardidaLio}), choose $\tilde{x}\in \Omega^+$ such that $\mathrm{dist}(\tilde{x},\Omega_0)<\mathrm{dist}(\tilde{x},\partial\Omega)$ and consider the ball $B_R=B_R(\tilde{x})\subset \Omega^+$ such that $\partial B_R (\tilde{x})\cap \partial\Omega_0\neq \emptyset$.

Observe that $f$ is a strictly increasing function on the interval $(0,\delta)$, $\delta<1$.

Fix a $x_0\in\partial B_R (\tilde{x})\cap \partial\Omega_0$, so $u(x_0)=0$ and $u(x)>0$ in $B_R=B_R(\tilde{x})$. Note that, up to diminishing $R$, we can suppose also that $u<\delta$ in $B_R(\tilde{x})$. Indeed, since $u(x_0)=0$ there exists a ball $B_{r_0}(x_0)$ such that $u<\delta$ in this ball; so by taking $R_1<R$ with $B_{R_1}(\tilde{x}_1)\subset B_{r_0}(x_0)$ for a point $\tilde{x_1}\in \Omega^+$, now just replace $\tilde{x},R$ by $\tilde{x}_1,R_1$.

Consider the annulus $E_R=B_R \setminus \overline{B}_{R/2}$ centered in $\tilde{x}$ and set $\mu :=\min_{\partial B_{R/2}} u \in (0,\delta)$.
We need to find a good barrier in $E_R$ for our nonlinear problem. This is accomplished in the following claim.

\begin{claim} \label{claim v radial vazquez}
There exists a nonnegative  classical subsolution $v\in C^2(\overline{E}_R)$ of
\begin{align*}
\left\{
\begin{array}{rclcc}
\mathcal{L}_1^- [v] &>& f(v)  &\mathrm{in} & \overline{E}_R \\
v &=& \mu  &\mbox{on} &\partial B_{R/2} \\
v &=& 0  &\mbox{on} &\partial B_{R}
\end{array}
\right.
\end{align*}
which is radially decreasing in $r=|x-\tilde{x}|$, convex and such that $\partial_\nu v >0$ on $\partial B_R$.
\end{claim}

\begin{proof}
We start choosing a large $\alpha >1$ such that
\begin{align}\label{escolha de alpha vazquez particular case}
\frac{1}{R^2}\,\{\, \alpha \,[\lambda (\alpha +1) +(n-1)\Lambda -\gamma R\,]-dR^2\, \}>C_0 \,\alpha \,,
\end{align}
where $C_0:=a \,( |\mathrm{ln}\mu|+2|\mathrm{ln} R|+|\mathrm{ln} 2|+|\mathrm{ln}(R/2)\,|)+m_0\,$, for $m_0:=\max_{s\in [0,1]} f(s)>0$.

Let $\varepsilon >0$ be such that
$
\varepsilon =\varepsilon (\alpha):= \mu \,\{ (R/2)^{-\alpha} - R^{-\alpha} \}^{-1} ,
$
and set
\begin{align}
v(x):=\varepsilon\, \{|x-\tilde{x}|^{-\alpha}-R^{-\alpha}\}.
\end{align}

\begin{rmk}{$($See, for instance, lemma 2.2 in \cite{arma2010}$.)$} \label{lema2.2arma}
Suppose that $u\in C^2(B)$ is a radial function, say $u(x)=\varphi (|x-x_0|)$, defined in a ball $B\subset \rn$ centered on $x_0$. Then $$\mathrm{spec}(D^2u (x)\,)=\left\{ \frac{\varphi ' (|x-x_0|)}{|x-x_0|}, \ldots, \frac{\varphi ' (|x-x_0|)}{|x-x_0|}, \varphi '' (|x-x_0|)  \right\}.$$
\end{rmk}

With this choice of $\varepsilon$, of course $v= \mu$ on $\partial B_{R/2}$ and $v=0$ on $\partial B_R$. Notice that $v(x)=\varphi (r)$ with $r=|x-\tilde{x}|$, thus $\varphi '(r)=-2\alpha\varepsilon\, r^{-\alpha-1} <0$ and so
\begin{align*}
\partial_\nu v(x)=Dv(x)\cdot \vec{\nu}=-\varphi ' (R) \frac{x-\tilde{x}}{|x-\tilde{x}|} \cdot \frac{x-\tilde{x}}{|x-\tilde{x}|} = -\varphi ' (R) >0
\end{align*}
for every $x\in \partial B_R=\partial B_R(\tilde{x})$, where $\vec{\nu}=-\frac{x-\tilde{x}}{|x-\tilde{x}|}\,$ is the interior unit normal to the ball $B_R$. Further, $\varphi ''(r)=\alpha (\alpha+1)\varepsilon \, r^{-\alpha-2}>0$, and by remark \ref{lema2.2arma} we have, in $\overline{E}_R$,
\begin{align*}
\mathcal{L}_1^-[v]&=\mathcal{M}_{\lambda,\Lambda}^-(D^2 v)-\gamma\, |Dv|-d\,v(x) \\
&= \alpha\varepsilon |x-\tilde{x}|^{-\alpha-2} \,\{ \lambda (\alpha +1) +(n-1)\Lambda\} -\gamma\alpha\varepsilon |x-\tilde{x}|^{-\alpha-1} - d\varepsilon |x-\tilde{x}|^{-\alpha}+d\varepsilon R^{-\alpha} \\
&\geq \varepsilon|x-\tilde{x}|^{-\alpha} \, R^{-2}\,\{\, \alpha \,[\lambda (\alpha +1) +(n-1)\Lambda -\gamma R\,]-dR^2\, \} \\
&> \alpha\, \varepsilon \,C_0\, |x-\tilde{x}|^{-\alpha}
\end{align*}
by the choice of $\alpha$ in \eqref{escolha de alpha vazquez particular case}. Now we claim that
\begin{align}\label{af f(v) vazquez particular}
f(v)\leq \alpha\, \varepsilon \,C_0\, |x-\tilde{x}|^{-\alpha} \;\; \mathrm{in}\;\;\overline{E}_R
\end{align}
and this will finish the proof of the claim \ref{claim v radial vazquez}. For $r=R$ this is obvious. Note that, for $r\neq R$, \eqref{af f(v) vazquez particular} is equivalent to
\begin{align}\label{af f(v) vazquez particular 2}
a \left(1-\frac{r^\alpha}{R^\alpha} \right) \left| \,\mathrm{ln}\varepsilon - \mathrm{ln} (r^{\alpha}) +  \mathrm{ln} \left(1-\frac{r^\alpha}{R^\alpha} \right) \right| \leq \alpha\, C_0.
\end{align}
But the left hand side of \eqref{af f(v) vazquez particular 2} is less or equal than
\begin{align*}
a\,|\mathrm{ln}\varepsilon|&+\alpha \,a\, |\mathrm{ln}r|+a y\, |\mathrm{ln}y|\\
&\leq a\,\{ \,|\mathrm{ln} \mu|+\alpha \,|\mathrm{ln} R|+|\mathrm{ln}(2^\alpha -1)|\,\} + \alpha a\, \{\, |\mathrm{ln}R|+ |\mathrm{ln}(R/2)|\, \}+a y |\mathrm{ln}y|\\
&\leq  \alpha\, C_0
\end{align*}
by using  $\alpha >1$ and definitions of $\varepsilon$, $m_0$, where $y:=1-\left( \frac{r}{R} \right)^\alpha \in \left[ 0,1-\frac{1}{2^\alpha} \right]$ for $r\in \left[\frac{R}{2},R\right]$.

\qedhere{\,\textit{Claim \ref{claim v radial vazquez}.} }
\end{proof}

By construction, $u\geq v$ on $\partial E_R=\partial B_R \cup \partial B_{R/2}$.

\begin{claim} \label{u geq v vazquez particular}
$u\geq v$ in $E_R$.
\end{claim}

\begin{proof}
Suppose not, i.e. that the open set $\mathcal{O}:= E_R\cap \{u<v\}$ is not empty.

Notice that since $v$ is decreasing in $r$ and $v=\mu <\delta$ on $r=R/2$, then both $u,v \in (0,\delta)$ in $E_R=B_R\setminus\overline{B}_{R/2}$, and we can use the monotonicity of $f$ in $(0,\delta)$ to obtain
\begin{align*}
\left\{
\begin{array}{rclcl}
\mathcal{L}_1^-[v]-\mathcal{L}_1^-[u] &>& f(v)-f(u) \geq  0 &\;\mbox{in}& \;\,\mathcal{O}\\
v &\leq & u &\;\mbox{on}& \partial\mathcal{O}\subset\partial E_R \cup \{v=u\}
\end{array}
\right.
\end{align*}
in the $C$-viscosity sense, and so the comparison principle (for example proposition 3.3 in \cite{Koike}), gives us that $v\leq u$ in $\mathcal{O}$, which contradicts the definition of $\mathcal{O}$.

\qedhere{\,\textit{Claim \ref{u geq v vazquez particular}.}}
\end{proof}

To finish the proof of lemma \ref{vazquez particular}, if $u\in C^1 (\Omega)$, we can just use the fact that $u$ has a minimum at the interior point $x_0$, so by claim \ref{u geq v vazquez particular} we have
$0=\partial_\nu u (x_0)\geq \partial_\nu v (x_0)>0$, a contradiction.

Assume  $u$ is only in $ LSC(\Omega)$. Let $\rho < R/2$.
Observe that $v$ is a $C^2$ function in $B_\rho(x_0)$ with $v<0\leq u$ in $B_\rho (x_0)\setminus E_R$ and, since $v\leq u$ in $E_R$ with $u(x_0)=v(x_0)=0$, then $v$ touches $u$ from below at the interior point $x_0$. Thus, by definition of $u$ being a $C$-viscosity supersolution on $\Omega$, we get
$$
0<\mathcal{L}^-[v](x_0)-dv(x_0)-f(v(x_0))=\mathcal{L}^-[v](x_0)=\mathcal{L}^-[v](x_0)-du(x_0)-f(u(x_0))\leq 0.
$$

\subsection{Proof of Lemma \ref{lemma2.3arma} }\label{appendix}

The proof follows the original idea of \cite{arma2010} for viscosity solutions, with the slight improvements that can be found in theorem 6.9 in \cite{Koike}.

\begin{proof} Inequalities \eqref{arma2010eq.uev} and \eqref{arma2010eq.uew} follow by a simple computation, by using the fact that $\mathrm{spec}(\xi\otimes \zeta)=\{0, \ldots ,0, \xi\cdot\zeta\}$, where $\xi\otimes \zeta \in \mathbf{M}_{n\times n}(\real)$, $(\xi\otimes \zeta)_{ij} :=\xi_i \,\zeta_j$ for all $\xi,\zeta \in \rn$.
Suppose, then, that $u$ is an $L^p$-viscosity solution of \eqref{arma2010eq.u}.

Let $\psi\in W^{2,p}_{\mathrm{loc}}(\Omega)$ such that $v-\psi$ attains a local maximum at $x_0$, namely $v\leq \psi$ and $v(x_0)=\psi (x_0)$, and let $\varepsilon>0$. Take some $\Omega '$ with $x_0\in \Omega '\subset\subset \Omega$ and set $a=\|u\|_{L^\infty (\Omega ')}$.

Define $\varphi =\frac{1}{m} \{ \mathrm{ln}(1+m\psi)\}\in W^{2,p}_{\mathrm{loc}}(\Omega)$, i.e. $\psi=\frac{1}{m}\{ e^{m\varphi} -1\}$ then, since $u-\varphi$ has a maximum at $x_0$, by using the definition of $u$ being an $L^p$-viscosity subsolution of \eqref{arma2010eq.u} and also \eqref{arma2010eq.uev} for the pair $\varphi,\psi$, we get
\begin{align*}
\frac{\mathcal{M}^+(D^2\psi)}{1+m\psi}\geq \mathcal{M}^+(D^2\varphi) +\mu |Du|^2 \geq f(x)-b(x)|D\varphi| -c(x)u-\tilde{\varepsilon}\quad\textrm{a.e. in } \mathcal{O}
\end{align*}
where $\tilde{\varepsilon}=\frac{\varepsilon}{e^{am}+1}$ and $\mathcal{O}\subset \Omega '$ is some subset with positive measure.

Let $\delta\in (0,1)$.
Notice that, since $(\psi -u) (x_0)=0$ and $\varphi \in W^{2,p}(\mathcal{O})\subset C(\overline{\mathcal{O}})$ for $p\geq n>n/2$,  there exists $\mathcal{O}_\delta\subset \mathcal{O}$ such that $v\leq \psi\leq v +\frac{\delta}{m}$ in $\mathcal{O}_\delta$. Thus $1+m\psi \leq 1+mv+\delta=e^{mu}+\delta$ and
\begin{align*}
\mathcal{M}^+(D^2\psi)+b(x)&|D\psi|\geq f^+(x)e^{mu} -f^-(x)(e^{mu}+\delta)\\
&+(c^+u^-+c^-u^+)e^{mu}-(c^+u^++c^-u^-)(e^{mu}+\delta) -{\varepsilon}   \quad\textrm{a.e. in } \mathcal{O}_\delta
\end{align*}
i.e. we have shown that $v$ is an $L^p$-viscosity subsolution of
\begin{align*}
\mathcal{M}^+(D^2v)+b(x)|D\psi|\geq f(x)e^{mu}-c(x)ue^{mu} -(f^- +c^+u^++c^-u^-)\delta   \quad\textrm{in } {\Omega}
\end{align*}
for any $\delta \in (0,1)$. The desired conclusion follows by letting $\delta\rightarrow 0$, since $\|f_\delta\|_{L^p(B)}\rightarrow 0$ for any $B\subset\subset\Omega$, where $f_\delta:=(f^- +c^+u^++c^-u^-)\delta$, by proposition \ref{Lpquad}. The proof of the remaining inequalities in the $L^p$-viscosity sense are similar.
\end{proof}

A direct consequence of lemma \ref{lemma2.3arma} is the following extension of proposition 3.5 in \cite{CCKS}.

\begin{corol} \label{2x diff a.e. quad.}
Let $f\in L^p_{\mathrm{loc}}(\Omega)$, $b\in L^\infty_+(\Omega)$, for $p\geq n$, $\mu\ge0$, and $u$ a locally bounded $L^p$-viscosity solution of
\begin{align*}
\mathcal{L}^+[u]+\mu |Du|^2 \geq f(x) \;\;\;\textrm{in } \Omega
\qquad \left( \,\mathcal{L}^-[u]-\mu |Du|^2 \leq f(x) \;\;\;\textrm{in } \Omega\, \right)
\end{align*}
then $u$ is twice superdifferentiable $($subdifferentiable$)$ a.e. in $\Omega$.
\end{corol}

\begin{proof}
Set  $v=\frac{1}{m}(e^{mu}-1)$ for $m=\frac{\mu}{\lambda_P}$. Say $b(x)\leq \gamma$. Then, by lemma \ref{lemma2.3arma} and $u\in L^\infty_{\mathrm{loc}} (\Omega)$, $v$ is an $L^p$-viscosity solution of
\begin{align*}
\mathcal{M}^+(D^2v)+\gamma |Dv|\geq f(x)(1+mv) \in L^p_{\mathrm{loc}}(\Omega).
\end{align*}
By proposition 3.5 in \cite{CCKS}, $v$ is twice superdifferentiable a.e. in $\Omega$, and so is $u$.
\end{proof}

As a matter of fact, corollary \ref{2x diff a.e. quad.} is true for $p>n-\varepsilon_0>n/2$ as in \cite{CCKS}, since lemma \ref{lemma2.3arma} also holds in this case.
Moreover, it follows from the argument in \cite{CCKS}, that locally bounded $L^p$-viscosity solutions of $F[u]+\langle M(x)D u,D u\rangle=f$, with $F$ satisfying the structure condition in \ref{SC0}\,, are twice differentiable a.e.  and  satisfy the equation at almost all points.

\end{document}